\documentclass[a4paper,11pt]{article}
\pdfoutput=1 

\usepackage{url}
\usepackage{amsmath}
\usepackage{amssymb}
\usepackage{amsfonts}
\usepackage{amssymb}
\usepackage{stmaryrd}
\usepackage{amscd}
\usepackage{epsfig}
\usepackage{epsf}
\usepackage{graphicx}
\usepackage{float}
\usepackage{dsfont}
\usepackage[ruled]{algorithm2e}
\usepackage{placeins}
\usepackage[latin1]{inputenc}
\usepackage[english]{babel}
\usepackage{theorem}
\usepackage[all]{xy}
\usepackage{multicol}
\usepackage{subfig}

\graphicspath{{Figures/}}

\def \mod {\mathop{\rm mod}\nolimits}
\def \p {{\lfloor n\alpha\rfloor}}

\def\eref#1{(\ref{#1})}

\def \N {\mathbb{N}}
\def \P {\mathbb{P}}
\def \A {\mathcal{A}}
\def \B {\mathcal{B}}
\def \C {\mathcal{C}}
\def \M {\mathcal{M}}
\def \V {\mathcal{V}}
\def \Z {\mathbb{Z}}
\def \z {\mathbf{0}}
\def \u {\mathbf{1}}
\def \p {\mathbf{?}}
\def \M {\mathcal{M}}
\def \A {\mathcal{A}}
\def \C {\mathcal{C}}
\def \N {\mathbb N}

\def \Z {\mathbb Z}
\def \P {{\mathbb P}}
\def \a {\alpha}

\def \env {\mbox{env}}

\newtheorem{defi}{Definition}[section]
\newtheorem{prop}[defi]{Proposition}
\newtheorem{lemm}[defi]{Lemma}
\newtheorem{cor}[defi]{Corollary}
\newtheorem{theo}[defi]{Theorem}

\newtheorem{conj}[defi]{Conjecture}
\theorembodyfont{\normalfont}
\newtheorem{exam}[defi]{Example}

\begin{document}

\title{Probabilistic cellular automata, invariant measures, \\ and
  perfect sampling\thanks{This work was partially supported by the ANR project MAGNUM (ANR-2010-BLAN-0204).}}

\author{Ana {\sc Bu\v{s}i\'c}%
\thanks{INRIA/ENS, 23, avenue d'Italie, CS 81321, 75214 Paris Cedex 13, France. 
E-mail: {\tt Ana.Busic@inria.fr}.}
\and
Jean {\sc Mairesse}%
\thanks{LIAFA, CNRS and Universit\'e Paris Diderot - Paris 7, Case 7014, 75205 Paris Cedex 13, France.
E-mail: {\tt Jean.Mairesse@liafa.jussieu.fr}.}
\and 
Ir\`ene {\sc Marcovici}%
\thanks{ENS Lyon/LIAFA, CNRS and Universit\'e Paris Diderot - Paris 7, Case 7014, 75205 Paris Cedex 13, France. E-mail: {\tt Irene.Marcovici@liafa.jussieu.fr}. }}

\date{\today}

\maketitle

\begin{abstract}

A probabilistic cellular automaton (PCA) can be viewed as a Markov chain. 
The cells are updated synchronously and independently, according to a distribution depending on a finite neighborhood.
We investigate the ergodicity of this Markov chain.
A classical cellular automaton is a particular case of PCA. For a
1-dimensional cellular automaton, we prove that ergodicity is equivalent to
nilpotency, and is therefore undecidable. We then propose an efficient perfect sampling 
algorithm for the invariant measure
of an ergodic PCA. Our algorithm does not assume any monotonicity
property of the local rule. It is based on a bounding process which
is shown to be also a PCA. Last, we focus on the PCA Majority, whose asymptotic behavior is unknown, and 
perform numerical experiments using the perfect
sampling procedure. 
\end{abstract}

\smallskip


\smallskip

{\noindent\bf AMS classification (2010):} Primary: 37B15, 60J05, 60J22. Secondary: 37A25, 60K35, 68Q80.
 

\tableofcontents


\section{Introduction}

A \emph{deterministic cellular automaton} (DCA) consists of a lattice (e.g. $\Z$ or $\Z^2$ or
$\Z/n\Z$) divided in regular cells, each cell containing a 
letter of a finite alphabet. The cells evolve
synchronously,  each one evolving in function of a
finite number of cells in its neighborhood, according to a local
rule.

DCA form a natural mathematical object: by Hedlund's theorem~\cite{hedl},
the mappings realized by DCA are precisely the continuous
functions (for the product topology) commuting with the shift. 
They also constitute a powerful model of computation, in particular they can
``simulate'' any Turing machine. Last, due to the amazing gap between
the simplicity of the definition and the intricacy of the generated behaviors,
DCA are good candidates for
modelling ``complex systems'' appearing in physical and biological
processes.

\medskip

\emph{Probabilistic cellular automata.} 
To take into account random events, one is led to consider 
 probabilistic versions of DCA. In one of them, 
at most one cell is updated at each time, this
cell being randomly chosen according to a given distribution. For an
infinite set of cells, one is led to consider continuous time
models, 
obtaining what is known in probability theory as an \emph{interacting
  particle system} \cite{ligg85}.  
Another model is that of \emph{probabilistic cellular automata}
(PCA)~\cite{toombook}. For PCA, time is discrete, 
and all the cells evolve synchronously as for DCA, but the difference is
that for each cell, the new content  is randomly chosen, independently of the others, according to
a  distribution depending only on a finite neighborhood of the
cell. 

Let us mention a couple of motivations. First, the investigation of
fault-tolerant 
computational models was the motivation for the Russian school
to study PCA~\cite{toombook,gacs}. Second, PCA appear in 
combinatorial problems related to the enumeration of directed
animals~\cite{dhar83,BoMe98,LeMa}. 
Third, in the context of the classification of DCA (Wolfram's program),
robustness to random errors can be used as a discriminating
criterion~\cite{FRST,schab}. 
Recently, PCA also proved to be pertinent for the density
classification problem~\cite{densitypb}, that is, testing efficiently
if some sequence contains more occurrences of $0$ or $1$. Last, PCA
are used in statistical physics and in life sciences. They modelize
various phenomena, from the dynamical properties of the neural
tissue~\cite{neuro} to competition between species. 

\medskip

We focus our study on the equilibrium behavior of PCA. Observe that
a PCA may be viewed as a Markov chain over the state
space $\A^E$, where $\A$ is the alphabet and $E$ is
the set of cells.  So the equilibrium is studied via the invariant
measures of the Markov chain. Several questions are in
order. 

\medskip

\emph{Ergodicity.} A PCA is \emph{ergodic} if it has a unique and
attractive invariant measure. 
A challenging problem in this area is the \emph{positive rates
  conjecture}. A PCA  is
said to have \emph{positive rates} if for any neighborhood, the updated
content of a cell can be any letter with a
strictly positive probability. The
positive rates conjecture states that any one-dimensional ($E=\Z$)
model with positive rates is ergodic. 
G\'acs exhibited in 2001
a very large and complex counter-example in a paper of more than 200
pages~\cite{gacs} which was published with an 
introductory article by Gray~\cite{gray01}. 
This counter-example is far from being completely understood and
several questions remain. For instance, does the conjecture hold true
for small alphabets and neighborhoods? Even for alphabets and
neighborhoods of size 2, the question is not settled. 

\medskip

\emph{Performance evaluation.} 
The second natural question is whether the invariant measures can be
evaluated. A PCA with an alphabet and a neighborhood of
size 2 is determined by four parameters. If the parameters satisfy a
given polynomial equation, there exists an invariant measure with an
explicit product form~\cite[Chapter 16]{toombook}. Under another
polynomial condition, there exists an
invariant measure with an explicit Markovian form, see \cite[Chapter 16]{toombook} or \cite{BoMe98}. 
What happens for generic values of the parameters? Or for PCA with a
larger neighborhood or alphabet? 
When explicit computation is not possible, simulation becomes the
alternative. 
Simulating PCA is known to be a
challenging task, costly both in time and space. Also, 
configurations cannot be tracked down one by one (there is an infinite
number of them when $E=\Z$) and may
only be observed through some measured parameters. 
So the crucial point is whether some guarantees can be given upon the
results obtained from simulations. 

\medskip

The contributions of the present paper are as follows: 

{\em First on ergodicity.} We prove that the ergodicity of a DCA on $\Z$ is
undecidable. This was mentioned as {\em Unsolved
  Problem 4.5} in \cite{toom01}. Since a DCA is a special case of a
PCA, it also
provides a new proof of the undecidability of the
ergodicity of a PCA (Kurdyumov, see \cite[Chap. 14]{toombook}, and
Toom~\cite{toom00}).

{\em Second on performance evaluation.} 
Given an ergodic PCA, a \emph{perfect sampling} procedure is a random algorithm which
returns a configuration distributed according to the invariant
measure. By applying the procedure repeatedly, we can
estimate the invariant measure with arbitrary precision. 
We propose such an algorithm for
PCA by adapting  the \emph{coupling from the past} method of
Propp $\&$ Wilson~\cite{PrWi}. When the set of cells is
$E=\Z/n\Z$, a PCA is a finite state space Markov chain. Therefore,
coupling from the past from all possible initial configurations provides a
basic perfect sampling procedure. But a very inefficient one since the
number of configurations is exponential in $n$. Here, the contribution 
consists in simplifying the procedure. We define a new PCA on an extended alphabet, called the
\emph{envelope PCA} (EPCA). We obtain a
perfect sampling procedure for the original PCA by running the EPCA on
a single initial configuration. When the set of cells is $E=\Z^d$, a
PCA is a Markov chain on an uncountable state space. So there is no
basic perfect sampling procedure anymore. We prove the following: If
the PCA is ergodic, then the EPCA may or may not be ergodic. If it is
ergodic, then we can use the EPCA to design an efficient perfect
sampling procedure (the result of the algorithm is the finite restriction of
a configuration with the right invariant distribution). In the case
$E=\Z$, we give a sufficient condition for the EPCA to be ergodic. 
The EPCA can be viewed as a systematic treatment of ideas
already used by Toom for \emph{percolation PCA} (see for
instance \cite[Section 2]{toom01}). 

The perfect sampling procedure can also be run on a PCA whose
ergodicity is unknown, with the purpose of testing it. 
We illustrate this approach on \emph{Majority}, prototype of a
PCA whose equilibrium behavior is not well understood. 
More precisely, we define a parametrized family of PCA, called 
\emph{Majority}$(\alpha)$, $\alpha \in (0,1)$. We prove that for
$\alpha$ large enough, the PCA has several invariant measures.  We conjecture the
existence of a phase transition  between two
situations: ({\it i}) several invariant measures; ({\it ii}) a unique but non-attractive
invariant measure. We provide some numerical
evidence for the phase transition, which would be the first example of
this kind. 
In fact, the mere 
existence of a PCA satisfying ({\it ii}) had been a
long standing open question  which was recently answered by the positive~\cite{ChMa10}.  

\medskip

Section \ref{se-pca} gives the basic definitions. Section \ref{sec:DCA} is devoted to the
ergodicity problem. Section \ref{sec:sampling} presents the perfect sampling
procedures. Last, Section \ref{se-majority} is devoted to the case study of the
Majority PCA. 

A short version without proofs of the paper appears in the
proceedings of the STACS'2011 Conference~\cite{BMMa}. 

\section{Probabilistic cellular automata}\label{se-pca}

%
Let $\A$ be a finite set called the \emph{alphabet}, and let $E$ be a
countable or finite set of \emph{cells}. We denote by $X$ the set $\A^E$ of
\emph{configurations}.  

We assume that $E$ is
equipped with a commutative semigroup structure, whose law is denoted
by $+$. In examples, we consider mostly the cases $E=\Z$ or
$E=\Z/n\Z$. 
Given $K\subset E$ and $V\subset E$, we define 
$$V+K = \bigl\{ u+v \in E \mid u \in V, v \in K
\bigr\} \:.$$

A \emph{cylinder} is a subset of $X$ having the form $\{x\in X \mid
\forall k\in K, x_k=y_k\}$ for a given finite subset $K$ of $E$ and a
given element 
$(y_k)_{k\in K}\in \A^K$. When there is no possible confusion, we
shall denote briefly by $y_K$ the cylinder $\{x\in X \mid \forall k\in K,
x_k=y_k\}$.  
For a given finite subset $K$, we denote by $\C(K)$
the set of all cylinders of base $K$.  

\medskip

Let us equip $X=\A^E$ with the product topology, which can be
described as the topology generated by cylinders. We denote by
$\M(\A)$ the set of probability measures on $\A$ 
 and by $\M(X)$ the set of probability measures on $X$ for the
 $\sigma$-algebra generated by all cylinder sets, which corresponds to
 the Borelian $\sigma$-algebra. For $x\in X$, denote by $\delta_x$ the
 Dirac measure concentrated on the configuration $x$. 

\medskip


\begin{defi}\label{de-pca}
Given a finite set $V \subset E$, a \emph{transition function} of \emph{neighborhood} $V$
is a function $f: \A^{V} \rightarrow  \M(\A)$. 
The \emph{probabilistic cellular automaton} (PCA) $P$ of transition
function $f$ is the application 
\begin{align*}
P: \M(X) & \rightarrow  \M(X) \\
 \mu & \mapsto \mu P \:,
\end{align*}
defined on cylinders by:
\[
\mu P(y_K)=\sum_{x_{V+K}\in \C(V+K)}\mu(x_{V+K})\prod_{k\in K}f((x_{k+v})_{v\in V})(y_k) \:.
\]
\end{defi}


Let us look at how $P$ acts on a Dirac measure $\delta_z$. 
The content $z_k$ of the $k$-th cell is
changed into the letter $a\in \A$ with probability $f((z_{k+v})_{v\in V})(a)$, independently of the evolution of the other cells.
The real
number $f((z_{k+v})_{v\in V})(a)
\in [0,1]$ is thus to be thought as the conditional
probability that, after application of $P$, the $k$-th cell
will be in the state $a$ if, before its application, the neighborhood
of $k$ was in the state $(z_{k+v})_{v\in V}.$ 

%

\medskip

Let $u$ be the uniform measure on $[0,1]$. We define the product measure $\tau=\bigotimes_{i\in E}u$ on $[0,1]^E$.

\begin{defi}\label{de-update}
 An \emph{update function} of the probabilistic cellular automaton $P$ is a
deterministic function $\phi:\A^E\times [0,1]^E\rightarrow \A^E$ (the
function $\phi$ takes as argument a configuration and a sample in $[0,
  1]^E$, and returns a new configuration), satisfying for each $x\in
\A^E$, and each cylinder $y_K$, 
$$\tau(\{r\in [0,1]^E; \phi(x,r)\in y_K\})=\prod_{k\in K}f((x_{k+v})_{v\in V})(y_k).$$ 
\end{defi}

In practice, it is always possible to define an update function $\phi$ for which the
value of $\phi(x,r)_k$ only depends on $(x_{k+v})_{v\in V}$ 
and on $r_k$. For example, if the alphabet is $\A=\{a_1,\ldots,a_n\}$, one can set
\begin{equation}\label{eq-update}
\phi(x,r)_k=\left\{
\begin{array}{l}
a_1 \ \mbox{ if } \ 0\leq r_k< f((x_{k+v})_{v\in V}
)(a_1)\\
a_2 \ \mbox{ if } \ f((x_{k+v})_{v\in V}
)(a_1)\leq r_k< f((x_{k+v})_{v\in V}
)(\{a_1, a_2\})\\
\vdots\\
a_n \ \mbox{ if } \ f((x_{k+v})_{v\in V}
(\{a_1,a_2,\ldots, a_{n-1}\})\leq r_k \leq 1.\\
\end{array}
\right.
\end{equation}

\begin{figure}
  \centering       
  \includegraphics[scale=1]{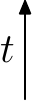}
    \; \, \; 
  \subfloat[$\alpha=0.5$]{\includegraphics[trim = 10cm 0.1cm 10cm 0.1cm, clip, scale=0.5]{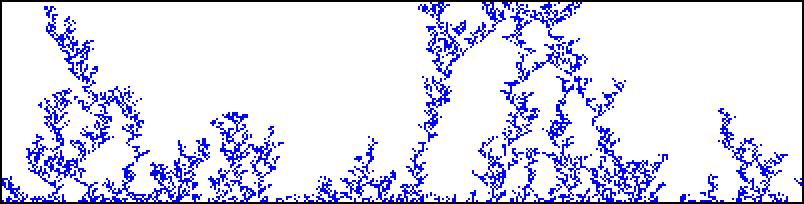}}        
  \; \, \;        
  \subfloat[$\alpha=0.6$]{\includegraphics[trim = 10cm 0.1cm 10cm 0.1cm, clip, scale=0.5]{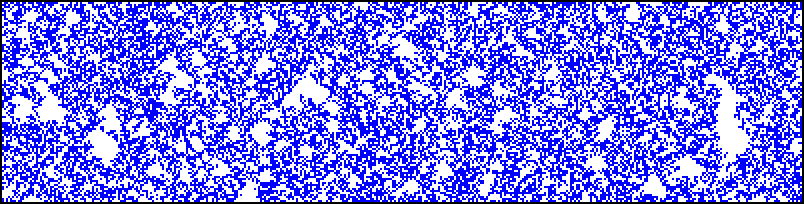}}
  \caption{\label{fi-pca1}
Space-time diagrams of the PCA of  Example \ref{xouy}, for $V=\{-1,0,1\}$
}
\end{figure}

For a given initial configuration $x^0\in\A^E$, and samples
$(r^t)_{t\in \N}, \ r^t \in [0,1]^E$, let $(x^t)_{t\in\N}\in {(\A^E)^\N}$
be the sequence defined recursively by
$$x^{t+1}=\phi(x^t, r^t).$$ 
Such a sequence is called a
\emph{space-time diagram}. It can be viewed as a realization of the
Markov chain. Examples of space-time diagrams appear in Figures \ref{fi-pca1} and
\ref{fi-exper}. 

\medskip
%
Classical cellular automata are a specialization of PCA.

\begin{defi}\label{de-dca}
A \emph{deterministic cellular automaton} (DCA) is a PCA such that for each sequence $(x_v)_{v\in V}\in \A^{V}$
, the measure $f((x_v)_{v\in V})$
is concentrated on a single letter of the alphabet. A DCA can thus be seen as a deterministic function 
$F:\A^E \rightarrow \A^E$.
\end{defi}

In the literature, the term \emph{cellular automaton}  denotes
what we call here a DCA.  Deterministic cellular automata have been widely studied, in particular
on the set of cells $E=\Z$, see Section~\ref{sec:DCA}. For a DCA, any initial configuration defines a unique space-time diagram.

\begin{exam}\label{somvois2} Let $\A=\{0,1\}$, $E=\Z$, and
  $V=\{0,1\}$. Consider $0<\varepsilon<1$ and the local function
\[
f(x,y)=(1-\varepsilon)\,\delta_{x+y \mod
    2}+\varepsilon\,\delta_{x+y+1 \mod 2} \:. 
\]
This defines a PCA that
can be considered as a perturbation of the DCA $F:\A^E\rightarrow \A^E$ defined by $F(x)_i=x_i+x_{i+1}
\mod 2$, with errors occurring in each cell independently with
probability $\varepsilon$. 
\end{exam}


\begin{exam}\label{xouy} 
Let $\A=\{0,1\}$, $E=\Z^d$, and let $V$ be a  finite subset of $E$. 
Consider $0<\a<1$ and the local function:
$$f((x_{v})_{v\in V})= \a\,\delta_{\max (x_v, \ v\in V)} +
  (1-\a)\,\delta_{0} \:.$$
The corresponding PCA is called the \emph{percolation PCA} associated with $V$ and
$\a$. The particular case of the space $E = \Z$ and the neighborhood $V = \{0,1\}$ 
is called the \emph{Stavskaya PCA}. In Figure \ref{fi-pca1}, we represent two space-time diagrams of the percolation PCA for $V=\{-1,0,1\}$.
\end{exam}

\subsection*{Invariant measures and ergodicity}\label{se-inva}

A PCA can be seen as a Markov chain on the state space $\A^E$. We use
the classical terminology for Markov chains that we now recall. 


\begin{defi}\label{de-statergo}
A probability measure $\pi\in\M(X)$ is said to be an \emph{invariant
  measure} of the PCA $P$ if $\pi
P=\pi$.  The PCA is \emph{ergodic} if it has exactly
one invariant measure $\pi$ which is \emph{attractive}, that is, for any measure $\mu\in\M(X)$, the
sequence $\mu P^n$ converges weakly to $\pi$ ({\it i.e.} for any cylinder $C$, $\lim_{n\rightarrow +\infty}\mu P^n(C)=\pi(C)$).
\end{defi}

A PCA has at least one invariant measure, and the set of invariant measures 
is convex and compact. This is a standard fact, based on the observation that the set $\M(X)$ of measures on $X$ is compact for the weak topology, see for instance \cite{toombook}. 
Therefore, there are three
possible situations for a PCA:

\begin{enumerate}
\item[({\it i})] several invariant measures;
\item[({\it ii})] a unique invariant measure which is not attractive;
\item[({\it iii})] a unique invariant measure which is attractive (ergodic case). 
\end{enumerate}

\begin{exam}\label{somvois2-2}
Consider the PCA of Example \ref{somvois2}. Using the results in
\cite[Chapters 16 and 17]{toombook}, one can prove that the PCA is
ergodic and that its unique invariant measure is the uniform mesure, {\it i.e.}
the product of Bernoulli measures of parameter 1/2. 
\end{exam}

\begin{exam}\label{xouy-2}
Consider the percolation PCA of Example \ref{xouy}. Observe that 
the Dirac measure $\delta_{0^{E}}$ is an invariant measure. 
Using a coupling with a percolation model, one can prove the
following, see for instance~\cite[Section 2]{toom01}. There exists $\alpha^* \in (0,1)$ such that:
\begin{eqnarray*}
\alpha < \alpha^* & \implies & (iii): \ \mbox{ergodicity} \\
\alpha > \alpha^* & \implies & (i): \ \mbox{several invariant
  measures.}
\end{eqnarray*}
The exact value of $\alpha^*$ is not known but it satisfies $1/|V| \leq \alpha^*
\leq 53/54$. 
\end{exam}

The existence of a PCA corresponding to situation ({\it ii}) had been a long
standing conjecture (
{\em Unsolved problem 5.7} in Toom~\cite{toom01}). In \cite{ChMa10}, it is proved that 
situation ({\it ii}) occurs for the PCA on $\{0,1\}^{\Z}$ with neighborhood
 $V=\{0,1\}$, and local function $a$
defined by $a(00)(1)=1/2,\: a(01)(1)=0,\: a(10)(1)=1,\: a(11)(1)=1/2$. 

The PCA of Example \ref{xouy} exhibits a phase transition between the situations ({\it i}) and ({\it iii}). In Section \ref{se-majority}, we study a PCA that may have a phase transition between the situations ({\it ii}) and ({\it iii}). It would provide the first example of this type.

\section{Ergodicity of DCA}
\label{sec:DCA}

DCA form the simplest class of PCA, it is therefore natural to study
the ergodicity of DCA. In this section, we prove the undecidability of ergodicity
for DCA (Theorem \ref{th-undec}). This also gives a new proof of the undecidability
of the ergodicity for PCA. 

\medskip

{\bf Remark.} In the context of DCA, the terminology of Definition
\ref{de-statergo} might be confusing. Indeed a DCA $P$ can be viewed
in two different ways: 
\[
(i) \mbox{ a (degenerated) Markov chain;} \qquad (ii) \mbox{ a symbolic 
dynamical system.}
\]
In the dynamical system terminology, 
$P$ is {\em uniquely ergodic} if: $[\exists ! \mu, \ \mu P=\mu].$
In the Markov chain terminology
(that we adopt), $P$ is {\em ergodic} if: $[\exists ! \mu, \ \mu
  P=\mu]$ and $[\forall \nu, \ \nu P^n \stackrel{w}{\longrightarrow}
  \mu]$, where $\stackrel{w}{\longrightarrow}$ stands for the weak
convergence. Knowing if the unique ergodicity (of symbolic dynamics)
implies the ergodicity (of the Markovian theory) is an open question for DCA. 

\medskip

The {\em limit set} of $P$ is defined by 
$$LS = \bigcap_{n\in \N} P^n({\A}^{E}) \:.$$
In words, a configuration belongs to $LS$ if it may occur after an
arbitrarily long evolution of the cellular automaton. 

Observe that $LS$ is non-empty since it is the decreasing limit of
non-empty closed sets. A constructive way to show that $LS$ is non-empty is as
follows. The image by $P$ of a monochromatic configuration $x^{E}$ is
monochromatic: $x^E \rightarrow y^E$. In particular there exists a
monochromatic periodic orbit for $P$, and we have:
\begin{equation}\label{eq-obs}
x_0^{E} \rightarrow x_1^{E} \rightarrow \cdots \rightarrow
x_{k-1}^E \rightarrow x_0^{E} \quad \implies \quad \{x_0^{E}, x_1^{E}, \dots , x_{k-1}^E \} \subset
LS \:.
\end{equation}

Recall that $\delta_u$ denotes the probability measure
concentrated on the configuration $u$. The periodic orbit $(x_0^{E},
\ldots , x_{k-1}^E)$ provides an invariant measure
given by $(\delta_{x_0^E}+\ldots+\delta_{x_{k-1}^E})/k$. More
generally, the support of any invariant measure is included in the
limit set. 

\begin{defi}\label{de-nil}
A DCA is {\em nilpotent} if its limit set is a singleton. 
\end{defi}

Using (\ref{eq-obs}), we see that a DCA is nilpotent iff
$LS=\{x^E\}$ for some $x\in \A$. The following stronger statement 
is proved in \cite{CPYu}, using a compactness argument: 
$$[ \ P \mbox{ \emph{nilpotent}} \ ] \quad \iff \quad [ \ \exists x \in \A,
  \exists N\in\N, \ \ P^N(\A^{E}) = \{x^E\} \ ] \:. $$

We get next proposition as a corollary. 

\begin{prop}
Consider a DCA $P$. We have:
\[
[ \ P \mbox{ nilpotent} \ ] \quad \implies \quad [ \ P \mbox{ ergodic}
  \ ] \:.
\]
\end{prop}

{\it Proof.} 
Let $x\in\A$ and $N\in\N$ be such that $P^N(\A^{E}) = \{x^E\}$. For any probability
measure $\mu$ on $\A^E$, we have $\mu P^N = \delta_{x^E}$. Therefore,
$P$ is ergodic with unique invariant measure $\delta_{x^E}$. 
$\hfill \square$ 

\medskip

If we restrict ourselves to DCA on $\Z$, we get the converse
statement. 

\begin{theo}\label{th-1d}
Consider a DCA $P$ on the set of cells $\Z$. We have:
\[
[ \ P \mbox{ nilpotent} \ ] \quad \iff \quad [ \ P \mbox{ ergodic}
  \ ] \:.
\]
\end{theo}

{\it Proof.} Let $P$ be an ergodic DCA. Assume that there exists a monochromatic
periodic orbit $(x_0^\Z, \ldots , x_{k-1}^\Z)$ with $k\geq 2$. Then 
$\mu= (\delta_{x_0^\Z} + \cdots + \delta_{x_{k-1}^\Z})/k$ is the
unique invariant measure. The sequence $\delta_{x_0^\Z} P^n$ does not
converge weakly to $\mu$, which is a contradiction. 
Therefore, there exists a monochromatic fixed point: $P(x^{\Z})=x^\Z$,
and $\delta_{x^\Z}$ is the unique invariant measure. 

Define the cylinder $C = \{v \in \A^\Z \mid \forall i \in K, \ v_i =
x\}$, where $K$ is some finite subset of $\Z$. For any initial
configuration $u \in \A^\Z$, using the ergodicity of $P$, we have:
\[
\delta_{u} P^n (C) \longrightarrow \delta_{x^\Z}(C) = 1 \:.
\]
But $\delta_{u} P^n$ is a Dirac measure, so $\delta_{u} P^n (C)$ is equal to $0$ or
$1$. Consequently, we have $\delta_{u} P^n (C)=1$ for $n$ large enough, that is,
$$\exists N\in\N, \forall n\geq N, \forall i \in K, \qquad P^n(u)_i = x \:.$$

In words, in any space-time diagram of $P$, any column becomes
eventually equal to $xxx\cdots$. Using the terminology of Guillon $\&$
Richard~\cite{GuRi}, the DCA $P$ has a {\em weakly nilpotent trace}. 
It is proved in \cite{GuRi} that the weak nilpotency of the trace
implies the nilpotency of the DCA. (The result is proved for 
cellular automata on $\Z$ and left open in larger dimensions.) 
This completes the proof. 
 $\hfill \square$ 

\medskip

Kari proved in \cite{kari} that the nilpotency of a DCA on $\Z$ is
undecidable. (For DCA on $\Z^d$, $d\geq 2$, the
proof appears in \cite{CPYu}.) By coupling Kari's result with Theorem
\ref{th-1d}, we get: 

\begin{cor}\label{th-undec}
Consider a DCA $P$ on the set of cells $\Z$. The ergodicity of $P$ in
undecidable. 
\end{cor}

The undecidability of the ergodicity of a PCA was a known
result, proved by 
Kurdyumov, see \cite{toombook}, see also Toom~\cite{toom00}.
Kurdyumov's and Toom's
proofs use a
non-deterministic PCA of dimension 1 and a reduction of the halting
problem of a Turing machine. 

Corollary \ref{th-undec} is a stronger statement. In fact, 
the (un)decidability of the ergodicity of a DCA was mentioned as {\em Unsolved
  Problem 4.5} in \cite{toom01}. 
We point out that Corollary \ref{th-undec} can also be obtained without Theorem \ref{th-1d}, by directly adapting Kari's proof
to show the undecidability of the ergodicity of the DCA associated with a NW-deterministic
tile set.







\section{Sampling the invariant measure of an ergodic PCA}
\label{sec:sampling}

Generally, the invariant measure(s) of a PCA cannot be described
explicitly. Numerical simulations are consequently very useful to 
get an idea of the
behavior of a PCA. Given an ergodic PCA, we propose a \emph{perfect
  sampling} algorithm which generates configurations
\emph{exactly} according to the invariant measure. 

\smallskip

A perfect sampling procedure for finite Markov chains has been
proposed by Propp $\&$ Wilson \cite{PrWi} using a \emph{coupling from the 
past} scheme. 
Perfect sampling procedures have been developed since in
various contexts. We mention below some works directly linked to the
present article. 
For more
  information see the annotated bibliography: \emph{Perfectly Random
    Sampling with Markov Chains},  
\url{http://dimacs.rutgers.edu/~dbwilson/exact.html/}.

\smallskip

The complexity of the algorithm depends on the number of all possible
initial conditions, which is prohibitive for PCA. 
A first crucial observation already
appears in \cite{PrWi}: 
for a monotone Markov chain, one has to consider 
two trajectories corresponding to minimal and maximal states of the
system. For anti-monotone systems, an analogous technique has been developed by H\"aggstr\"om $\&$ Nelander \cite{haggnel}
that also considers 
only extremal initial conditions.
To cope with more general
situations, Huber \cite{hube} introduced 
the idea of a bounding chain for determining when coupling has occurred.
The construction of these 
bounding chains is model-dependent and in general not straightforward.   
In the case of a Markov chain on a lattice, Bu\v{s}i\'c et
al. \cite{BGVi} proposed  
an algorithm to construct bounding chains. 

\smallskip

Our contribution is to show that the bounding chain ideas can be given in a
particularly simple and convenient form in the context of PCA
via the introduction of the \emph{envelope PCA}. 

\subsection{Basic coupling from the past for PCA}\label{se-basicalgo}

\subsubsection{Finite set of cells}

Consider an ergodic PCA $P$ on the alphabet $\A$ and on a finite set of cells $E$ (for example $\Z_m=\Z/m\Z$). Let $\pi$ be
the invariant measure on $X=\A^E$. A {\em perfect sampling} procedure is a random
algorithm which returns a state $x\in X$ with probability
$\pi(x)$. Let us present the Propp $\&$ Wilson, or {\em coupling from
  the past (CFTP)}, perfect sampling procedure. 

\begin{algorithm}[H]\label{algobf} 
\KwData{An update function $\phi:X\times [0,1]^{E} \rightarrow X$
  of a PCA. A family $(r^{-n}_k)_{(k,n)\in E \times \N}$ of i.i.d. r.v. with uniform distribution in $[0,1]$.}
\Begin{
$t=1$ \;
\Repeat{$|R_0|=1$}{
$R_{-t}=X$ \;
\For{$j=-t$ \KwTo $-1$}{$R_{j+1}=\{\phi(x,(r_i^{j})_{i\in E})\: ;\: x\in R_{j}\}$}
$t=t+1$}
\Return{the unique element of $R_0$}}
\caption{Basic CFTP algorithm for a finite set of cells}
\end{algorithm}

The good way to implement this algorithm is to keep track of the partial couplings of trajectories. This allows to consider only one-step transitions.

%

\begin{prop}[\cite{PrWi}]\label{pr-cftp}
If the procedure stops almost surely, then the PCA is ergodic and the
output is distributed according to the invariant measure. 
\end{prop}

\begin{figure}
\begin{center}
\includegraphics[scale=0.5]{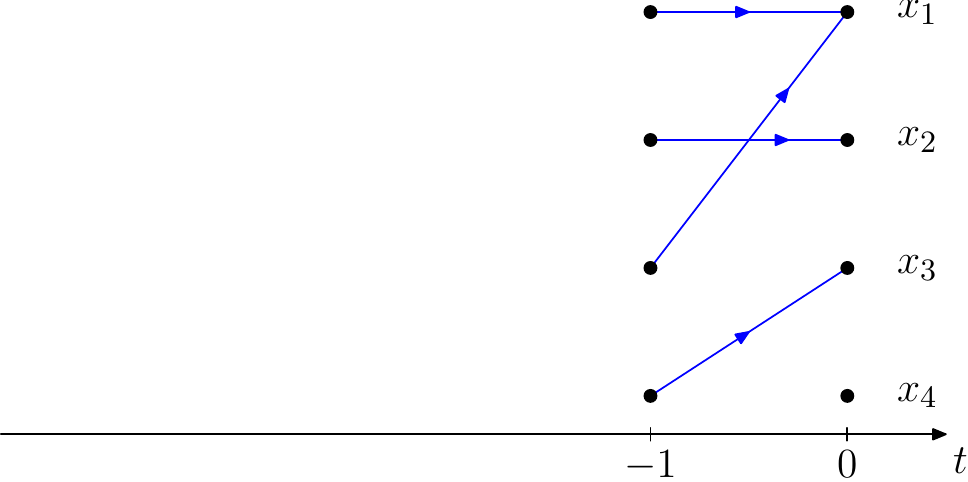} \qquad  \qquad 
\includegraphics[scale=0.5]{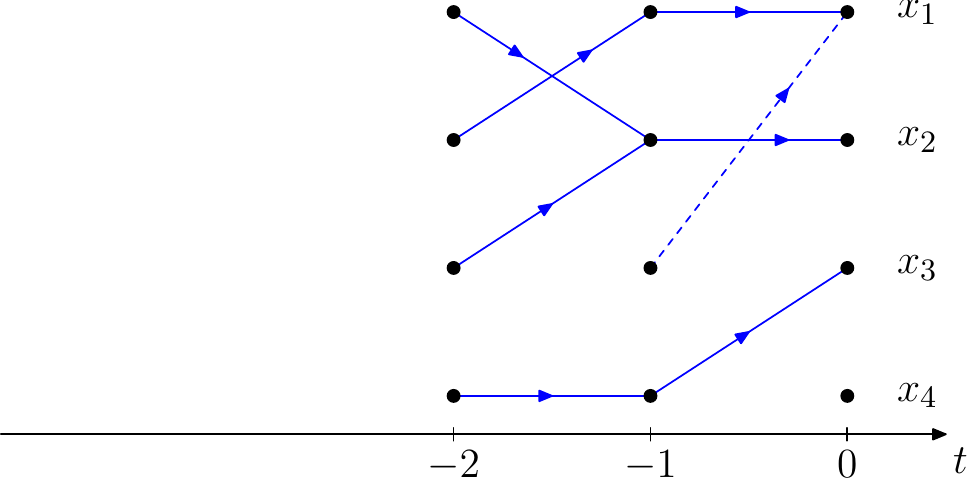}

\bigskip

\,

\bigskip

\includegraphics[scale=0.5]{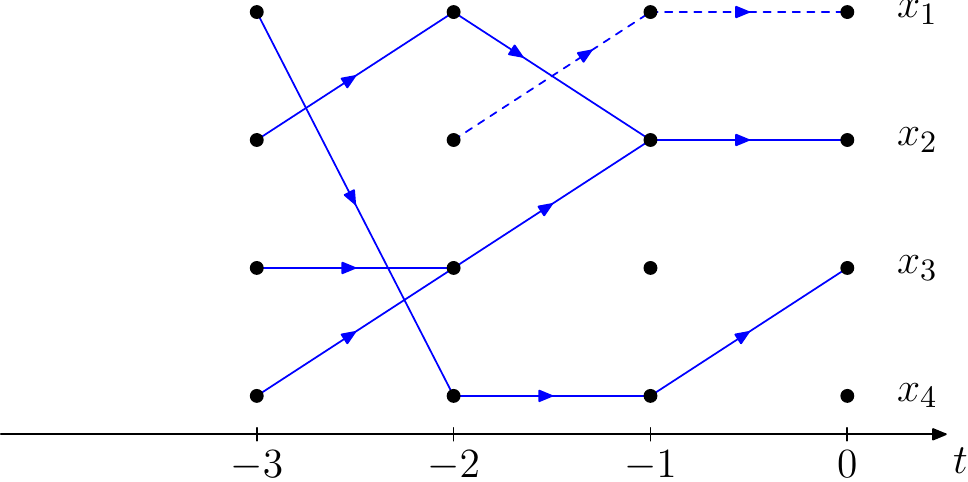}
\qquad \qquad 
\includegraphics[scale=0.5]{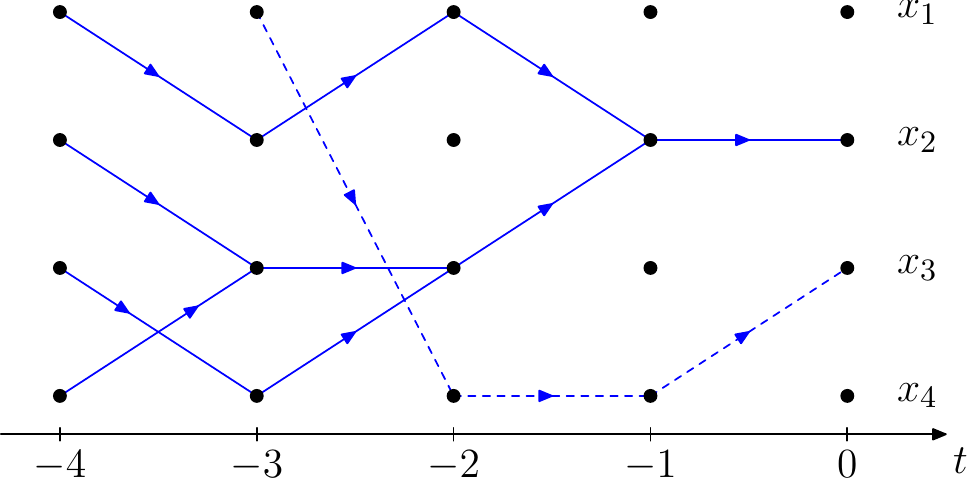}
\caption{Coupling from the past}\label{fi-cftp}
\end{center}
\end{figure}

In Figure \ref{fi-cftp}, we illustrate the algorithm on the toy example of a PCA on the
alphabet $\{0,1\}$ and the set of cells $\Z_2$. The state space is thus
$X=\{x_1=00, x_2=01,x_3=10, x_4 = 11\}$. On this sample, the algorithm
returns $x_2$.

A sketch of the proof of Proposition \ref{pr-cftp} can be given using 
Figure \ref{fi-cftp}. On the last of the four pictures, the
Markov chain is run from time -4 onwards and its value is $x_2$ at
time 0. If we had run the Markov 
chain from time $-\infty$ to 0, then the result would obviously still
be $x_2$. But if we started from time $-\infty$, then the Markov chain
would have reached equilibrium by time 0. 

\subsubsection{Infinite set of cells}

Assume that the set of cells $E$ is infinite. Then a PCA defines a Markov chain on the infinite state
space $X=\A^{E}$, so the above procedure is not effective anymore. 
However, it is possible to use the locality of the updating rule of a
PCA to still define a perfect sampling procedure. (This observation
already appears in \cite{BeSt99}.)

\medskip

Let $P$ be an ergodic PCA and denote by $\pi$ its invariant
distribution. In this context, a {\em perfect sampling} procedure is a
random algorithm taking as input a finite subset $K$ of $E$ and
returning a cylinder $x_K \in \C(K)$ with probability $\pi(x_K)$. 

\begin{figure}[H]
\begin{center}
\includegraphics[scale=0.7]{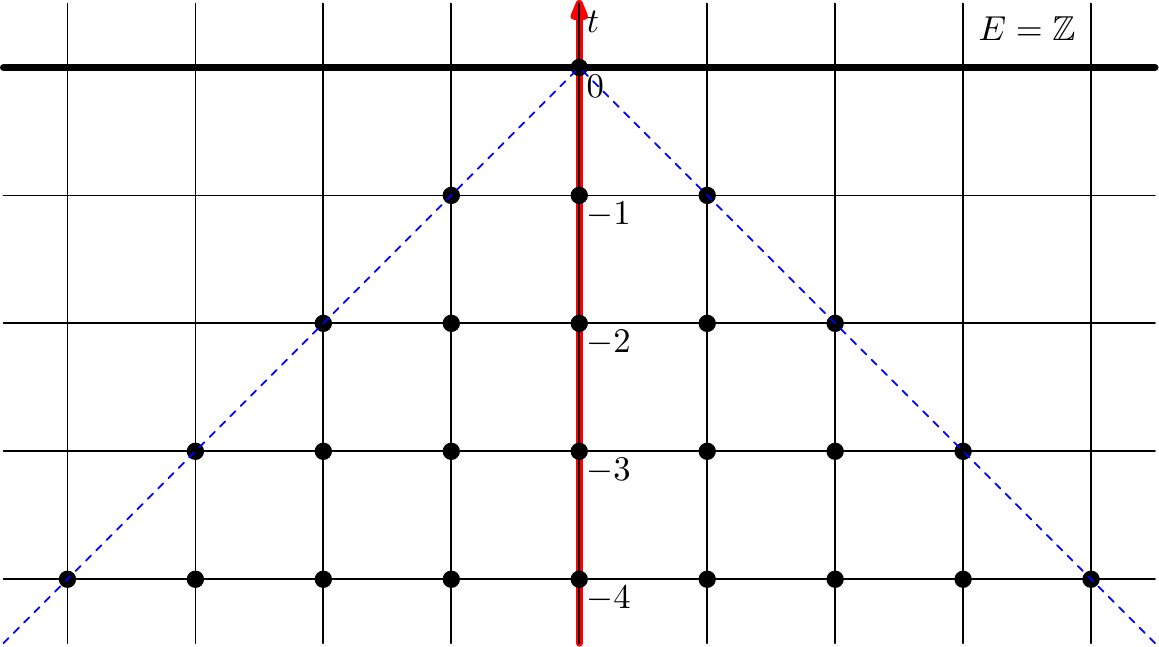}
\caption{Dependence cone of a cell}\label{fi-cone}
\end{center}
\end{figure}

To get such a procedure, we use the following fact: if the PCA is run from time $-k$
onwards, then to compute the content of the cells in $K$ at time 0, it is
enough to consider the cells in the finite dependence cone of $K$.
This is illustrated in Figure \ref{fi-cone} for the set of cells $E=\Z$ and the neighborhood $V=\{-1,0,1\}$, with the choice $K=\{0\}$.

Observe that the orientation has changed with respect to
Figure \ref{fi-cftp} in order to be consistent with the convention of Figures
\ref{fi-pca1} and \ref{fi-exper} for space-time diagrams. 

\medskip

Let us define this more formally. Let $V$ be the neighborhood of the PCA. Given a subset $K$ of $E$,
the \emph{dependence cone} of $K$ is the family $(V_{-t}(K))_{t\geq 0}$ of subsets of $E$ defined recursively by $V_0(K)=K$ and
$V_{-t}(K) =V+V_{-t+1}(K)$. Let $\phi:X\times [0,1]^E\rightarrow X$ be an update function, for instance the one defined
according to \eref{eq-update}. For a given subset $K$ of $E$, we denote $\phi_{-t}:\A^{V_{-t}(K)} \times
[0,1]^{V_{-t}(K)}\rightarrow \A^{V_{-t+1}(K)}$ the corresponding restriction of
$\phi$. 

With these  notations, the algorithm can be written as follows.


\begin{algorithm}[H]\label{algobi}
\KwData{An update function $\phi:X\times [0,1]^E\rightarrow X$ of a
  PCA. A
  family $(r^{-n}_k)_{(k,n)\in E\times \N}$ of i.i.d. r.v. with
  uniform distribution in $[0,1]$. A finite subset $K$ of $E$.}
\Begin{
$V_0(K)=K$ \;
$t=1$ \;
\Repeat{$|R_0|=1$}{
$V_{-t}(K)=V+V_{-t+1}(K)$ \;
$R_{-t}={\A}^{V_{-t}(K)}$ \;
\For{$j=-t$ \KwTo $-1$}{$R_{j+1}=\{\phi_j(x,(r^{j}_i)_{i\in V_j(K)})\: ;\: x\in R_{j}\}\subset \A^{V_{j+1}(K)}$}
$t=t+1$}
\Return{the unique element of $R_0$}}
\caption{Basic CFTP algorithm for an infinite set of cells}
\end{algorithm}

Next proposition is an easy extension of Proposition \ref{pr-cftp}.

\begin{prop}\label{pr-cftpinfty}
If the procedure stops almost surely, then the PCA is ergodic and the
output is distributed according to the marginal of the invariant
measure. 
\end{prop}

\subsection{Envelope probabilistic cellular automata (EPCA)}

The CFTP algorithm is 
inefficient when the state space is large. This is the case for
PCA: when $E$ is finite, the set $\A^E$ is very large, and when $E$ is
infinite, it is the number of configurations living in the dependence cone described above which is very
large. We cope with this difficulty by introducing the {\em envelope}
PCA. 

\medskip

To begin with, let us assume that $P$ is a PCA on the alphabet
$\A=\{0,1\}$ (as previously, the set of cells is denoted by $E$, the
neighborhood by $V\subset E$, and the local function by $f$). The case
of a general alphabet is treated in Section \ref{sse-exte}. 

\paragraph{Definition of the EPCA.}


Let us introduce a new alphabet:
$${\B}=\{{\z},{\u},{\p}\}.$$ 
A word on $\B$ is to be
thought as a word on  $\A$  in which the letters
corresponding to some positions are not known, and are thus replaced
by the symbol ``${\p}$''. Formally we identify $\B$ with
$2^{\A}-\emptyset$ as follows: ${\z} = \{0\}$, ${\u}= \{1\}$, and
${\p}=\{0,1\}$. So each letter of $\B$ is a set of possible letters of
$\A$. With this interpretation, we view
 a word on $\B$ as a set of words on $\A$. For instance, $$\p\u\p = \{
 010,011,110,111 \}.$$ 

We will associate to the PCA $P$ a new PCA on the alphabet $\B$, that we call the \emph{envelope probabilistic cellular automaton} of $P$.

\begin{defi}\label{de-epca} 
The \emph{envelope probabilistic cellular automaton (EPCA)} of $P$, is the PCA $\env(P)$ of alphabet $\B$, defined on the set of cells $E$, with the same neighborhood $V$ as for $P$, and a local function 
$\env(f):{\B}^V\rightarrow {\M}({\B})$
defined for each $y\in {\B}^V$ by

$$\env(f)(y)({\z})=\min_{x\in \A^V,\; x\in y}f(x)(0)$$
$$\env(f)(y)({\u})=\min_{x\in \A^V,\;  x\in y}f(x)(1)$$
$$\env(f)(y)({\p})=1-\min_{x\in \A^V, \; x\in y}f(x)(0)-\min_{x\in \A^V,\; x\in y}f(x)(1).$$
\end{defi}

We point out that $\min_{x\in \A^V,\; x\in y}f(x)(1)+\max_{x\in \A^V,\; x\in y}f(x)(0)=1$, so that the last quantity $\env(f)(y)({\p})$ is non-negative.

Moreover, $\env(P)$ acts like $P$ on configurations which do not
contain the letter ``${\p}$''. More precisely, 
\begin{equation}\label{eq-coincide}
\forall y\in {\A}^V, \qquad \env(f)(y)({\z})=f(y)(0), \quad
\env(f)(y)({\u})=f(y)(1), \quad \env(f)(y)({\p})=0 \:.
\end{equation}

In particular, we get the
following.

\begin{prop}\label{le-same}
If the EPCA $\env(P)$ is ergodic then the PCA $P$ is ergodic. 
\end{prop}

{\it Proof.} 
According to \eref{eq-coincide}, any invariant measure of $P$ corresponds to an
invariant measure of $\env(P)$. Therefore, if $P$ has several
invariant measures, so does $\env(P)$. Assume that $P$ has a unique invariant
measure $\mu$ which is non-ergodic. Let $\mu_0$ be such that
$\mu_0P^n$ does not converge to $\mu$. Then $\mu_0\ \env(P)^n$ does not
converge either, see \eref{eq-coincide}. To summarize, we have proved
that $P$ non-ergodic implies $\env(P)$ non-ergodic. 
$\hfill \square$ 

\medskip

The converse of Proposition \ref{le-same} is not true
 and counter-examples will be given in Section
\ref{sse-cex}. 

\paragraph{Construction of an update function for the EPCA.}

Let us define the update function
$$\tilde{\phi}:\B^E\times [0,1]^E \rightarrow \B^E$$
of the PCA $\env(P)$, by:
\begin{equation}\label{eq-update2}
\tilde{\phi}(y,r)_k=\left\{
\begin{array}{l}
{\z} \mbox{ if } 0\leq r_k< \env(f)((y_{k+v})_{v\in V})({\z})\\
{\u} \mbox{ if } 1-\env(f)((y_{k+v})_{v\in V})({\u})\leq r_k\leq1 \\
{\p} \mbox{ otherwise.}\\
\end{array}
\right.
\end{equation}

The value of $\tilde{\phi}(y,r)_k$ in function of $r_k$ can thus be represented as follows.
\begin{center}
\includegraphics[scale=0.8]{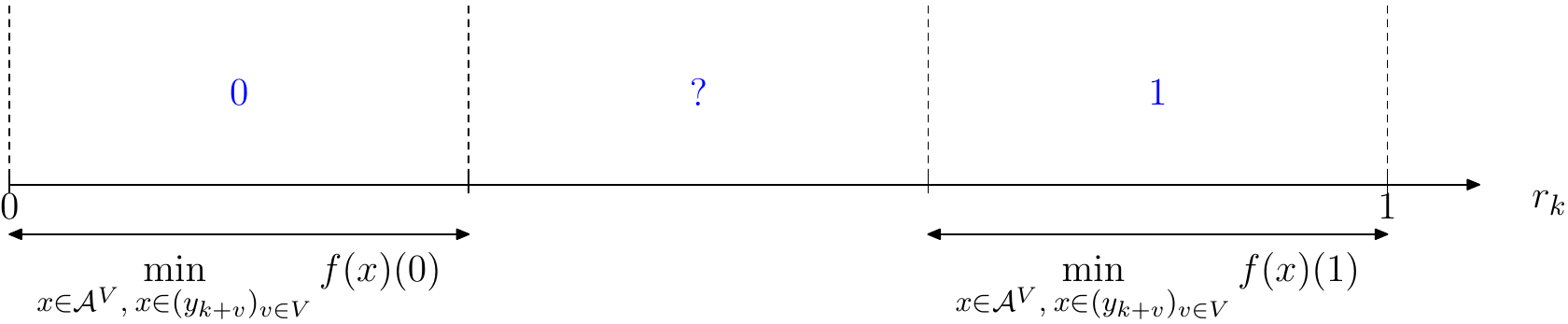}
\end{center}

For a PCA of neighborhood $V=\{0,1\}$, we represent below the construction of the updates of the EPCA when the value of the neighborhood is $\z ?$.

\begin{center}
\includegraphics[scale=0.65]{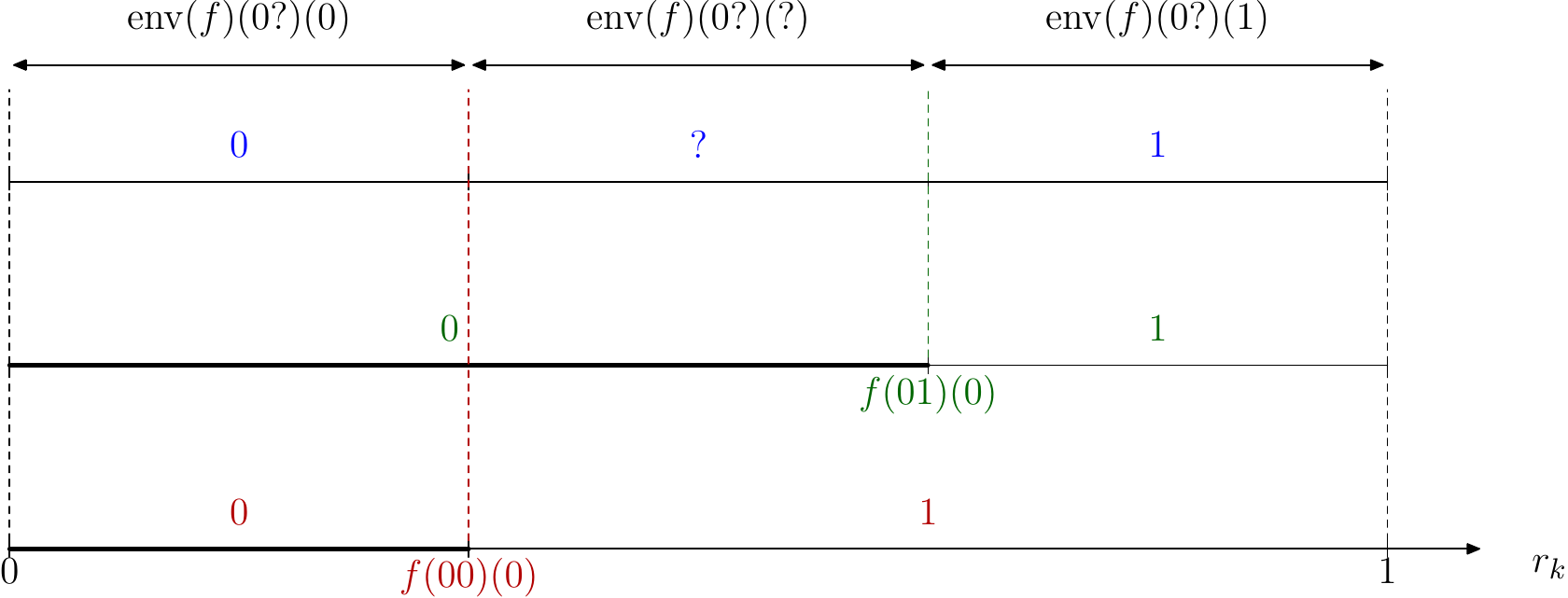}
\end{center}

\medskip


Let $\phi$ be the natural update function for the PCA $P$ defined 
as in \eref{eq-update}. 
Observe that  the function $\tilde{\phi}$ coincides with $\phi$ on
configurations which do not contain the letter
``${\p}$''. Furthermore, we have:
\begin{equation}\label{eq-updates}
\forall r\in [0,1]^E, \ \forall x\in \A^E, \ \forall y\in \B^E, \qquad
x\in y \implies \phi(x,r)\in \tilde{\phi}(y,r) \:.
\end{equation}

\subsection{Perfect sampling using EPCA}\label{se-sampl.epca}

We propose two perfect sampling algorithms, for a finite and for an infinite number of cells.
We show that in both cases, the algorithm stops almost surely if and
only if the EPCA is ergodic (Theorem \ref{prop:EPCAalgos}).
The ergodicity of the EPCA implies the ergodicity of the PCA but the converse is not true:
we provide a counterexample for each case, finite and infinite
(Section \ref{sse-cex}).
We also give conditions of ergodicity of the EPCA (Prop. \ref{pr-cvepcafini} and \ref{pr-suff}).

\subsubsection{Algorithms}\label{sse-algo}


\paragraph{Finite set of cells.} The idea is to consider only one trajectory of the EPCA - the one that starts from the initial configuration ${\p}^E$ (coding the set of all configurations of the PCA). 
The algorithm stops when at time $0$, this trajectory hits the set $\A^E$. 

\begin{algorithm}[H]\label{algoef}
\KwData{The pre-computed update function $\tilde{\phi}$. A family $(r^{-n}_k)_{(k,n)\in E\times \N}$ of i.i.d. r.v. with uniform distribution in $[0,1]$.}
\Begin{
$t=1$ \;
\Repeat{$c \in \A^E$}{
$c={\p}^E$ \;
\For{$j=-t$ \KwTo $-1$}{$c=\tilde{\phi}(c,(r^{j}_i)_{i\in E})$}
$t=t+1$}
\Return{$c$}}
\caption{Perfect sampling using the EPCA for a finite set of cells}
\end{algorithm}

\paragraph{Infinite set of cells.}

Once again, we consider only one trajectory of the EPCA.

\begin{algorithm}[H]\label{algoei}
\KwData{The pre-computed update function $\tilde{\phi}$. A family
  $(r^{-n}_k)_{(k,n)\in E\times \N}$ of i.i.d. r.v. with
  uniform distribution in $[0,1]$. A finite subset $K$ of $E$.}
\Begin{
$V_0(K)=K$ \;
$t=1$ \;
\Repeat{$c \in \A^K$}{
$V_{-t}(K)=V+V_{-t+1}(K)$ \;
$c={\p}^{V_{-t}(K)}$ \;
\For{$j=-t$ \KwTo $-1$}{$c=\tilde{\phi}_j(c,(r^{j}_{i})_{i\in V_j(K)})\in\B^{V_{j+1}(K)}$}
$t=t+1$}
\Return{$c$}}
\caption{Perfect sampling using the EPCA for an infinite set of cells}
\end{algorithm}

\begin{theo}
\label{prop:EPCAalgos}
Algorithm \ref{algoef}, resp. \ref{algoei}, stops almost surely if and only if the EPCA is ergodic. In that case, the output of the algorithm 
is distributed according to the unique invariant measure of the PCA. 
\end{theo}


{\it Proof.} The argument is the same in the finite and infinite
cases. We give it for the finite case. 
Assume first that Algorithm \ref{algoef} stops
almost surely. By construction, it implies that for all $\mu_0$, the
measure $\mu_0\ \env(P)^n$ is asymptotically supported by
$\A^E$. Therefore, we can strengthen the result in Proposition
\ref{le-same}: the invariant measures of $\env(P)$ coincide with the
invariant measures of $P$. In that case, $\env(P)$ is ergodic iff $P$ is
ergodic. 
Using \eref{eq-updates}, the halting of Algorithm \ref{algoef} implies the
halting of Algorithm \ref{algobf}. 
Furthermore, if we use the same samples $(r^{-n}_{k})_{(k,n)\in E\times \N}$,
Algorithms \ref{algoef} and \ref{algobf}  will have the same output. 
According to Proposition \ref{pr-cftp}, this output is
distributed according to the unique invariant measure of $P$. In
particular, $P$ is ergodic. So $\mbox{env}(P)$ is ergodic. 

\medskip

Assume now that the EPCA is ergodic. The unique invariant measure
$\pi$ of
$\env(P)$ has to be supported by $\A^E$. Also, by ergodicity, we have
$\delta_{\p^E} \ \env(P)^n \stackrel{w}{\longrightarrow} \pi$. This
means precisely that the Algorithm \ref{algoef} stops a.s. 
\hfill $\square$

\subsubsection{Criteria of ergodicity for the EPCA}

\paragraph{Finite set of cells.}

In the next proposition, we give a necessary and sufficient condition for
the EPCA to be ergodic.
In particular, this condition is satisfied if the PCA has positive rates (see the Introduction).

\begin{prop}\label{pr-cvepcafini}
The EPCA $\env(P)$ is ergodic  if and only if $\env(f)({\p}^V)({\p})<1$. This condition can also be written as: 
\begin{equation}\label{eq-necsuff}
\min_{x\in \A^V} f(x)(0)+\min_{x\in\A^V} f(x)(1)>0.
\end{equation}
\end{prop}

{\it Proof.} If $\env(f)({\p}^V)({\p})=1$, then for almost any
$r\in[0,1]^E$, we have $\tilde{\phi}({\p}^E,r)={\p}^E,$ so that at
each step of the algorithm, the value of $c$ is ${\p}^E$ with
probability 1. 

Conversely, if we assume for example that $p=\min_{x\in
  {\A}^V}f(x)(0)>0$, then for any configuration $d\in \B^E$, the
probability to have $\tilde{\phi}(x,r)={\z}^E$ is greater than $p^{|E|}$,
so that the algorithm stops almost surely, and the expectation of the running time can be roughly bounded by $1/p^{|E|}$. \hfill $\square$

\paragraph{Infinite set of cells.}

For an infinite set of cells the situation is more complex. The
condition of Proposition \ref{pr-cvepcafini} is not sufficient to
ensure the ergodicity of the EPCA. A counter-example is
given in Section \ref{sse-cex}. 
First, we propose a rough sufficient condition of convergence for
Algorithm \ref{algoei}. 

\begin{prop}\label{pr-suff}
Let $\alpha^*\in (0,1)$ be the critical probability of the percolation PCA
with neighborhood $V$, see Examples \ref{xouy} and \ref{xouy-2}. 
The EPCA $\env(P)$ is ergodic if 
\begin{equation}\label{critergod}
\env(f)(\p^V)({\p}) \ < \ \alpha^*
\end{equation}
and non-ergodic if
\begin{equation}\label{critnonergod}
\min_{x\in \B^V-\A^V} \env(f)(x)({\p}) \ > \ \alpha^* .
\end{equation}
\end{prop}

\medskip

{\it Proof.} 
Recall that $\B = \{{\bf 0}, {\bf 1}, \p \}$. Define $\C= \left\{ {\bf d}, \p \right\}$, with ${\bf d} =
\{ {\bf 0}, {\bf 1} \}$. A word over $\C$ is interpreted as a
set of words over $\B$, for instance, ${\bf d}? = \{{\bf 0}\p, {\bf
  1}? \}$. The symbol {\bf d} stands for {\bf d}etermined letter, as
  opposed to ?  which represents an unknown letter. 

We define a new PCA $Q$ on the alphabet $\C$, with the same
neighborhood $V$ as $P$ and $\env(P)$, and with the transition
function $g:{\C}^V\rightarrow {\M}({\C})$ defined by:
\[
g({\bf d}^V) = \delta_{{\bf d}}, \quad \mbox{and} \quad \forall u \in \C^V - \{{\bf
  d}^V\}, \quad g(u) = \alpha 
\delta_{\p} + (1-\alpha) \delta_{{\bf d}} \:,
\]
for $\alpha = \max_{x\in \B^V} \env(f)(x)({\p}) =
\env(f)(\p^V)({\p})$. 

Observe that $\delta_{{\bf d}^E}$ is an invariant measure of $Q$. 
Recall that $\tilde{\phi}$ is an update function of $\env(P)$, see
\eref{eq-update2}. Given the way $Q$ is defined, 
we can construct an update function $\phi_Q$ of $Q$ such that
\begin{equation}\label{eq-link}
\forall x \in \B^E, \forall y \in \C^E, \forall r\in [0,1]^E, \qquad x
\in y \implies \tilde{\phi}(x,r) \in \phi_Q(y,r) \:.
\end{equation}
In particular, assume that $Q$ is ergodic. Then $\delta_{\p^E} \ Q^n
\stackrel{w}{\longrightarrow} \delta_{{\bf d}^E}$. Using
\eref{eq-link}, it implies that Algorithm \ref{algoei} stops almost
surely, and
$\env(P)$ is ergodic according to Theorem \ref{prop:EPCAalgos}. To summarize, the ergodicity of $Q$ implies the ergodicity
of $\env(P)$. 

Observe that the PCA $Q$ is a percolation PCA as defined in Example \ref{xouy} 
(here, {\bf d} plays the role of 0 and $\p$ plays the
role of 1). Let $\alpha^*\in (0,1)$ be the critical probability of the percolation PCA
with neighborhood $V$, see Example \ref{xouy-2}. For $\alpha <
\alpha^*$, the percolation PCA $Q$ is ergodic. This completes the
proof of \eref{critergod}. 

\medskip

Define a PCA $R$ on the alphabet $\C$, with 
neighborhood $V$, and with the transition
function:
\[
h({\bf d}^V) = \delta_{{\bf d}}, \quad \mbox{and} \quad \forall u \in \C^V - \{{\bf
  d}^V\}, \quad h(u) = \beta 
\delta_{\p} + (1-\beta) \delta_{{\bf d}} \:,
\]
for $\beta = \min_{x\in \B^V-\A^V} \env(f)(x)({\p})$. 
Given the way $R$ is defined, 
we can construct an update function $\phi_R$ of $R$ such that
$$\forall x \in \B^E, \forall y \in \C^E, \forall r\in [0,1]^E, \forall k\in E, \qquad
        [ x \in y, \ \phi_R(y,r)_k=\p ] \implies \tilde{\phi}(x,r)_k = \p
        \:. 
$$
Therefore, the ergodicity of $env(P)$ implies the ergodicity of
$R$. Equivalently, the non-ergodicity of $R$ implies the
non-ergodicity of $env(P)$. 
Observe that the PCA $R$ is a percolation PCA. Therefore, for $\beta >
\alpha^*$, the percolation PCA $R$ is non-ergodic. This completes the
proof of \eref{critnonergod}.

\hfill $\square$

\subsubsection{Counter-examples}\label{sse-cex}

Recall Proposition \ref{le-same}: [ EPCA ergodic ] $\implies$ [ PCA ergodic
]. We now show that the converse is not true. 

\medskip

\begin{exam}
Consider the PCA with alphabet $\A=\{0,1\}$, neighborhood $V=\{-1,0,1\}$,
set of cells $E=\Z/n\Z$, and 
transition function 
\[
f(x,y,z) = \begin{cases} \delta_{1-y} & \mbox{if } xyz \in \{101,010\} \\
                         \alpha \delta_y + (1-\alpha) \delta_{1-y} &
                         \mbox{otherwise} \:,
\end{cases}
\]
for a parameter $\alpha \in (0,1)$. This is the PCA Majority studied in Section \ref{se-majority}.
For $n$ odd, we prove in Proposition \ref{pr-evenodd} that the PCA is
ergodic. However the associated EPCA satisfies $\env(f)(???) =
\delta_{?}$. According to Proposition \ref{pr-cvepcafini}, the EPCA is not
ergodic. 
\end{exam}

%

\begin{exam}
Consider the PCA of Example \ref{somvois2}. This PCA has positive
rates, in particular, it satisfies 
\eref{eq-necsuff}. So the EPCA is ergodic on a finite set of
  cells. Now let the set of cells be $\Z$. 
The PCA is ergodic for $\varepsilon \in (0,1)$, see Example
\ref{somvois2-2}. Consider the associated EPCA $\env(P)$. 
Assume for instance that $\varepsilon \in (0,1/2)$. We have 
\[
\env(f)(u) = \begin{cases} f(u) & \mbox{if } u \in \{\z,\u\}^V \\
\varepsilon \delta_{{\bf 0}} + \varepsilon \delta_{{\bf 1}} +
(1-2\varepsilon) \delta_{\p} & \mbox{otherwise}\:.
\end{cases}
\]
By applying Proposition \ref{pr-suff}, $\env(P)$ is non-ergodic if $1-2\varepsilon > \alpha^*.$
\end{exam}

\subsection{Decay of correlations}\label{sse-shift}

In what follows, the set of cells is $E=\Z^d, \ d\geq 1$. It is easy
to prove that the invariant measure of an ergodic PCA is
shift-invariant. Using the coupling from the past tool, we give
conditions for the invariant measure of an ergodic PCA to be
shift-mixing. 

\begin{defi} A measure $\mu$ on $X= {\cal A}^{\Z^d}$ is \emph{shift-mixing} if for
  any non-trivial translation shift $\tau$ of $\Z^d$, and for any cylinders $U,V$ of $X$,
\begin{equation}\label{shiftm}
\lim_{n\rightarrow +\infty}\mu(U\cap \tau^{-n}(V))=\mu(U)\mu(V).
\end{equation}
\end{defi}

The proof of the following
proposition is inspired from the proof of the validity of the coupling
from the past method (see \cite{PrWi} or \cite{haggb}). 

\begin{figure}[H]
\begin{center}
\includegraphics[scale=.5]{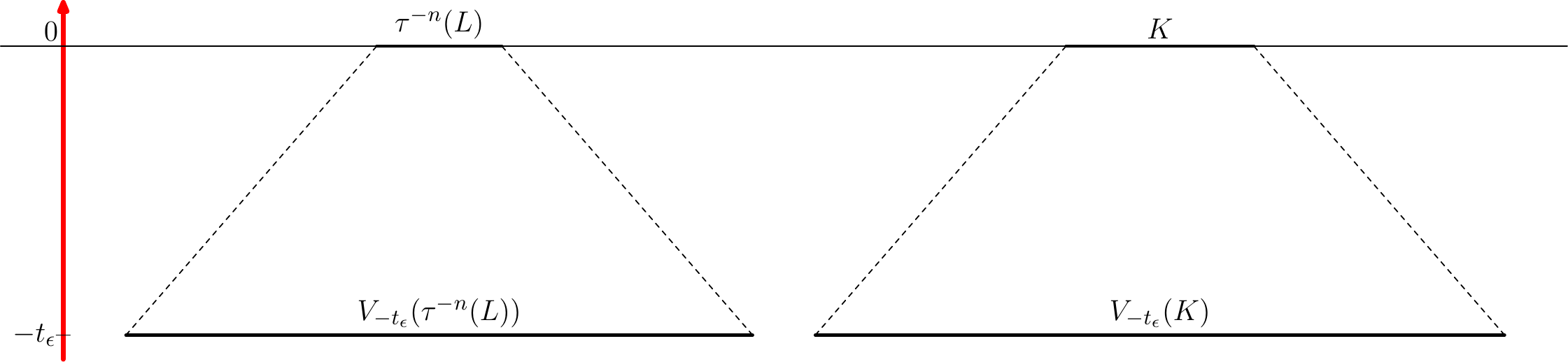}
\caption{Illustration of the proof of Proposition \ref{pr-decay}}\label{fi-shift}
\end{center}
\end{figure}

\begin{prop}\label{pr-decay}
If Algorithm \ref{algobi} stops almost surely, then the unique invariant measure of the PCA is shift-mixing. It is in particular the case under condition (\ref{critergod}).
\end{prop}

{\it Proof.} Assume that $P$ is an ergodic PCA, and denote by $\pi$
its unique invariant measure. Let $K$ and $L$ be two finite subsets of
$E$, and denote by $x_K$ and $y_L$ some cylinders corresponding to
these subsets. Since the perfect sampling algorithm stops almost
surely, for each $\varepsilon>0$, there exists an integer $t_\varepsilon$
such that with probability greater than $1-\varepsilon$, the algorithm
stops before reaching the time $-t_{\varepsilon}$ when it is run for the
set of cells $K$ or for the set of cells $L$. If $n\in \N^d$ is large enough, 
we have: $V_
{-t_{\varepsilon}}(K)\cap V_
{-t_{\varepsilon}}(\tau^{-n}(L))=\varnothing$.  

Let $Z$ be the output of the algorithm if it is asked to sample the
marginals of $\pi$ corresponding to the cells of $K\cup
\tau^{-n}(L)$. 

Imagine running the PCA from time $-t_{\varepsilon}$ and set of cells
$V_ {-t_{\varepsilon}}(K)\cup V_ {-t_{\varepsilon}}(\tau^{-n}(L))$
up to time $0$, using the same update variables as the ones used to get $Z$. Choose the initial condition at time
$-t_{\varepsilon}$ as follows: independently on  $V_
{-t_{\varepsilon}}(K)$ and $V_ {-t_{\varepsilon}}(\tau^{-n}(L))$,
and according to the relevant marginals of $\pi$.  
Let $X$, resp. $Y$, be the output at time 0 on the set of cells $K$,
resp. $\tau^{-n}(L)$. Observe that $X$ and $Y$ are distributed
according to the marginals of $\pi$. Furthermore, $X$ and $Y$ are
independent since the dependence cones of $K$ and $\tau^{-n}(L)$
originating at time $-t_{\varepsilon}$ are
disjoint. 
We therefore get
\begin{align*}
\pi(x_K\cap \tau^{-n}(y_L))-\pi(x_K)\pi(y_L)&= \P(Z_K=x_K, Z_{\tau^{-n}(L)}=y _L)-\P(X=x_K) \P(Y=y _L)\\
&= \P(Z_K=x_K, Z_{\tau^{-n}(L)}=y _L)-\P(X=x_K,Y=y _L)\\
&\leq  \P((Z_K, Z_{\tau^{-n}(L)})=(x_K, y _L) \mbox{ and } (X,Y)\not=(x_K,y _L))\\
&\leq  \P((Z_K, Z_{\tau^{-n}(L)})\not= (X,Y)) \ 
\leq \ 2\varepsilon \:.
\end{align*}
In the same way, we get $\pi(x_K)\pi(y_L)-\pi(x_K\cap
\tau^{-n}(y_L))\leq 2\varepsilon$. It completes the proof. \hfill
$\square$ 

\medskip

In Proposition \ref{pr-decay}, the coupling from
the past method is not used as a sampling tool but as a way to
get theoretical results. 
Knowing if there exists an ergodic PCA having an invariant measure which is not shift-mixing is an open question (see \cite{ColTiss10} for details).

\subsection{Extensions}\label{sse-exte}

In a PCA, the dynamic is homogeneous
in space. It is possible to get rid of this characteristic by defining
non-homogeneous PCA, for which the neighborhood and the transition function depend on the position of the cell. The definition below is to be compared with Definition \ref{de-pca}. The configuration space $X = \A^{E}$ is unchanged. 

%

\begin{defi} For each $k\in E$, denote by $V_k\subset E$ the (finite) neighborhood of the cell $k$, and by $f_k : \A^{V_k} \rightarrow \M(\A)$ the transition function associated to $k$. Set $\V(K)=\cup_{k\in K} V_k$. The {\em non-homogeneous PCA (NH-PCA)} of transition functions $(f_k)_{k\in E}$ is the application $P:
\M(X)\rightarrow \M(X), \ \mu \mapsto \mu P$, defined on cylinders by
\[
\mu P(y_K)=\sum_{x_{\V(K)}\in \C(\V(K))}\mu(x_{\V(K)})\prod_{k\in K} f_k((x_v)_{v\in V_k}
)(y_k)\:.
\]
\end{defi}


Observe that it is not necessary for $E$ to be equipped with a
semigroup structure anymore. We use this below to define the
finite restriction of a PCA. 

\medskip

It is quite straightforward to adapt the coupling from the past
algorithms to NH-PCA. More precisely, given a
NH-PCA, we define the associated NH-EPCA by
considering Definition \ref{de-epca} and replacing $V$ and $\env(f)$ by $V_k$ and 
$\env(f)_k$ for each $k\in E$. 
The algorithms of Section \ref{se-basicalgo} and \ref{sse-algo} are then unchanged, and
Proposition \ref{le-same} and Theorem \ref{prop:EPCAalgos} still hold in the non-homogeneous setting. 

\medskip

In Section \ref{se-majority}, we study the PCA Majority by approximating it by a
sequence of NH-PCA. Let us explain the construction in
a general setting. 

\medskip

Let $P$ be a PCA on the infinite set of cells $E$, with neighborhood $V$
and transition function $f:\A^V\rightarrow \M(\A)$. Let $D$ be a
finite subset of $E$. Define
\[
\overline{V}(D)= \{ u \in E \mid \exists x \in D, \exists v \in V\cup \{0\}, \ u=x+v\},
\quad B(D) = \overline{V}(D) - D \:.
\]
The set $B(D)$ is the \emph{boundary} of the domain $D$. Fix a probability
measure $\nu$ on $\A$. The {\em restriction} of $P$ associated with
$\nu$ and $D$ is the NH-PCA $P(\nu,D)$ with set of cells $\overline{V}(D)$ and
neighborhoods:
\[
\forall u \in D, \ V_u = \{u\} + V, \qquad \forall u \in B(D), \ V_u = \varnothing \:;
\]
and transition functions:
\[
\forall u \in D, \ f_u = f, \qquad \forall u \in B(D), 
\ f_u(\cdot) = \nu\:.
\]
In words, the boundary cells are i.i.d. of law $\nu$ and
the cells of $D$ are updated according to $P$. 

\medskip

If $\mu$ is a probability measure on $\A^S$, where $S$ is a finite
subset of $E$, we define 
its extension $\tilde{\mu}$ on $\A^E$ by setting, for a fixed letter $a\in \A$: 
\[
\forall x\in \A^E, \ \tilde{\mu}(x)=
\begin{cases}
\mu((x_k)_{k\in S}) & \mbox{if } \forall i\in E - S, \ x_i=a\\
0 & \mbox{  otherwise.}
\end{cases}
\]

\begin{lemm}\label{pr-approx}
Let $(D_i)_{i\in \N}$ be an increasing sequence of finite domains
$D_i\subset E$ such that $\cup_{i\in\N}D_i=E$. Let $(\nu_i)_{i\in \N}$ be a
sequence of probability measures on $\A$. For each $i$, let $\mu_i$
be an invariant measure of $P(\nu_i,D_i)$. 
Any accumulation point of the sequence
  $(\tilde{\mu}_i)_{i\in\N}$ is an invariant measure of the original
  PCA $P$ defined on $E$. 
\end{lemm}

{\it Proof.} 
Upon extracting a subsequence, we may assume that
$(\tilde{\mu}_j)_{j\in\N}$ converges to 
$\tilde{\mu}\in\M(X)$. 
We need to prove that for any cylinder $y_K\in\C(K),$ we
have $\tilde{\mu} P(y_K)=\tilde{\mu}(y_K)$.

By definition, $\mu_j P(\nu_j,D_j) =\mu_j$. 
Let the subset $K$ of $E$ and the cylinder $y_K\in\C(K)$ be fixed. If 
$j$ is large enough, we have $K\subset D_j$ and $\overline{V}(K)\subset
D_j$. So that $\mu_j(y_K)=\tilde{\mu}_j(y_K)$ and $P(\nu_j,D_j)$ and
$P$ coincide on $K$. We deduce that 
$\tilde{\mu}_j P(y_K) =\tilde{\mu}_j(y_K)$.  
By
taking the limit on both sides, we get
$\tilde{\mu} P(y_K)=\tilde{\mu}(y_K)$. 
$\hfill \square$ 

\subsubsection*{Alphabet with more than two elements}

The EPCA and the associated algorithms have been defined on a two
letters alphabet. It is possible to extend the approach to a general
finite alphabet. 

\medskip

Let $\A$ be the finite alphabet. Let $P$ be a PCA
with set of cells $E$, neighborhood $V$, and transition function $f:
\A^V \rightarrow \M(\A)$.  

Consider the alphabet $\B = 2^\A - \{\varnothing\}$, that is, the set of
non-empty subsets of $\A$. A word over $\B$ is viewed as a set of
words over $\A$.

\medskip

The EPCA $\env(P)$ associated with $P$ is the PCA on the alphabet $\B$
with neighborhood $V$ and transition function $\env(f)$ defined by: 
\[
\forall v \in \B^V, \ \forall y \in \B, \qquad \env(f)(v)(y) = \sum_{x
  \subset y} (-1)^{|y|-|x|} \min_{u \in v} f(u)(x) \:.
\]
For instance, we have:
$\env(f)(v)(\{0,1,2\}) = \alpha_{0,1,2} - \alpha_{0,1} - \alpha_{1,2} -
\alpha_{0,2} + \alpha_0 + \alpha_1 + \alpha_2$
with $\alpha_S = \min_{u \in v} f(u)(\{S\})$. 

\medskip

The algorithms of Section \ref{se-sampl.epca} are unchanged. Observe however that
the construction of an update function is not as natural as in 
the two-letters alphabet case. 

\section{The majority PCA: a case study}\label{se-majority}

The \emph{Majority} PCA is one of the simplest examples of PCA whose behaviour is not well understood. Therefore, it provides a good case study
for the sampling algorithms of Section \ref{sec:sampling}.

\subsection{Definition of the majority PCA}

\begin{figure}[H]
\begin{center}
\includegraphics[scale=0.8]{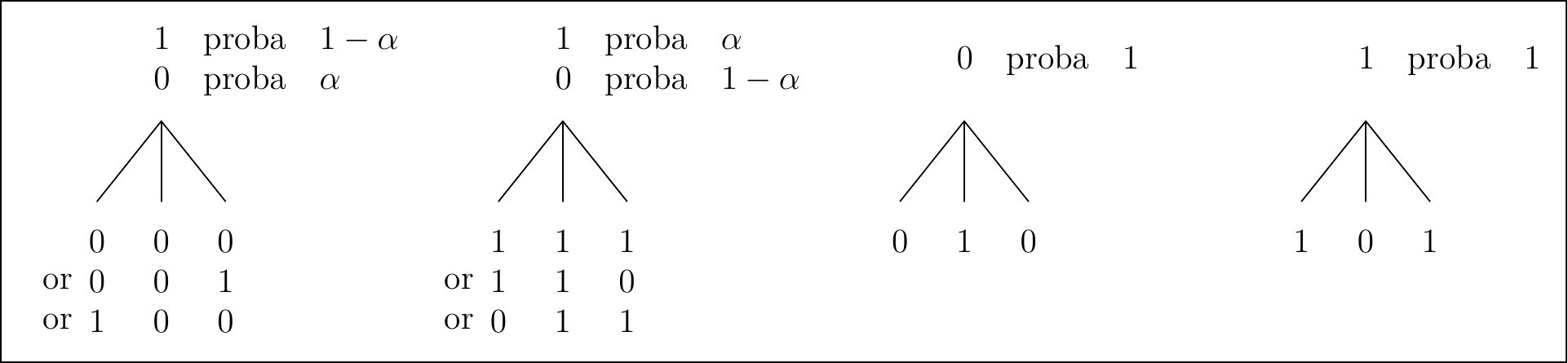}
\caption{The transition function of the PCA Majority}\label{maj}
\end{center}
\end{figure}

Given $0<\a<1$, the PCA \emph{Majority}$(\alpha)$, or simply {\em Majority}, 
is the PCA on the alphabet $\A=\{0,1\}$, with set of cells $E=\Z$ (or $\Z/n\Z$),
neighborhood  $V=\{-1,0,1\}$, and transition function
$$f(x,y,z)= \a\,\delta_{\mbox{maj}(x,y,z)}+ (1-\a)\,\delta_{1-y}\:,$$
where $\mbox{maj}:\A^3 \rightarrow \A$ is 
 the \emph{majority function}: the value of $\mbox{maj}(x,y,z)$ is
  $0$, resp. 1,  if
  there are two or three $0$'s, resp $1$'s, in the sequence $x,y,z$. The transition function of PCA Majority$(\alpha)$ can thus be represented as in Figure \ref{maj}. It consists in choosing
independently for each cell to apply the elementary rule 232 (with probability $\alpha$) or to flip the value of the cell.

\medskip

The PCA Minority$(\alpha)$ has also been studied (see \cite{schab}). It is defined by the transition function $g(x,y,z)=f(1-x,1-y,1-z)$.

\medskip

Let $x=(01)^{\Z}\in \{0,1\}^\Z$ be 
defined by:  $\forall n\in \Z$, $x_{2n}=0,
\ x_{2n+1}=1$.  The configuration  $(10)^{\Z}$ is defined
similarly. Consider the probability measure
\begin{equation}\label{eq-trivial}
\mu= (\delta_{(01)^{\Z}} +
\delta_{(10)^{\Z}})/2\:.
\end{equation}
Clearly $\mu$ is an invariant measure for the PCA Majority. 
The question is whether there exists other invariant measures. 

%
%
%
%



To get some other insight on the question of invariant measures, consider the PCA 
Majority on the set of cells $\Z_n={\Z / n\Z}$. This PCA 
has two completely
different behaviors depending on the parity of
$n$. 

\begin{prop}\label{pr-evenodd} Consider the Markov chain on the state space
  $\{0,1\}^{\Z_n}$ which is induced by the Majority PCA. The  Markov chain has a unique invariant measure $\nu$. If 
$n$ is even then $\nu=(\delta_{(01)^{n/2}} +\delta_{(10)^{n/2}})/2$; if
$n$ is odd then $\nu$ is supported by $\{0,1\}^{\Z_n}$. 
\end{prop} 

{\it Proof.} Let $Q$ be the transition matrix of the Markov chain.
We are going to prove the following. If $n$ is odd, then $Q$ is irreducible and aperiodic. 
If $n$ is even, the graph of $Q$ has a unique terminal component 
consisting of the two points $(01)^{n/2}=0101\ldots 01$ and
$(10)^{n/2}=1010\ldots 10$. 

\medskip

Assume first that $n$ is odd. From the configurations
$0^{n}=0\ldots 0$ or $1^{n}=1\ldots 1$, we can go to any configuration in one step.
So, we only need to prove that from a given
configuration, say $c\in \{0,1\}^{\Z_n}$, it is possible to reach
$0^{n}$ or $1^{n}$.
We assume that $c$ is neither $0^{n}$ nor $1^{n}$, and we
consider the length $L$ of the longest subsequence of identical 
bits in $c$. Since $n$ is odd, there are at least two consecutive
$0$'s or two consecutive $1$'s. So, we have $L\geq 2$. 
Say that we have the pattern $10^L1$ in $c$. With probability
$(1-\a)^2\a^L$, the update of
the corresponding cells will be: flip, majority, majority, ...,
majority, flip. Thus after one step, we can obtain configurations with the pattern
$0^{L+2}$ (or $0^{n}$ if $L=n-1$). By iterating, we see that after at most $n/2$ steps, we can reach
configuration $0^{n}$ or $1^{n}$. 

Assume now that $n$ is even. Clearly 
\[
\delta_{(01)^{n/2}}Q =
\delta_{(10)^{n/2}}, \qquad \delta_{(10)^{n/2}}Q =
\delta_{(01)^{n/2}}\:.
\]
Thus, $\{(01)^{n/2},(10)^{n/2}\}$ is a
terminal component of the graph of $Q$. Consider
$c \in \{0,1\}^{\Z_n} - \{(01)^{n/2}, (10)^{n/2}\}$. 
Then $c$ has at least two consecutive $0$'s or $1$'s.
We can use the same argument as for $n$ odd, to show that 
we can reach $0^{n}$ or $1^{n}$ in at most $n/2$ steps. 
We conclude by observing that $(01)^{n/2}$ can be reached from
$0^{n}$ or $1^{n}$ in one step. $\hfill \square$ 

\medskip

Let us come back to the PCA Majority on $\Z$. The invariant measure
$\mu$ in \eqref{eq-trivial} can be viewed as the ``limit'' over $n$
of the invariant measures of the PCA on $\Z_{2n}$. What about the
``limits'' of the invariant measures of the PCA on $\Z_{2n+1}$~? Do
they define other invariant measures for the PCA on $\Z$~? 

\medskip

One of the motivations of our work on perfect sampling algorithms for PCA was to test the following conjecture, which is inspired by the observations made in \cite{regn08} and \cite{schab} on a PCA equi\-valent to Majority. This conjecture concerns the existence of a ``phase transition'' phenomenon for the PCA Majority. 

\begin{conj}\label{co-co} 
There exists $\a_c \in (0,1)$ such that Majority$(\a)$ has a unique
invariant measure if $\a<\a_c$, and several invariant measures if
$\a>\a_c$.
\end{conj}

In the next subsection, we give some rigorous (but partial) results about the invariant measures of Majority($\a$). We first introduce a related PCA and use it to prove that if $\alpha$ is large enough, Majority$(\a)$ has indeed non-trivial invariant measures; we then present a dual model that could be used to have some information for small values of $\alpha$.

The last subsection is devoted to the experimental study of Majority$(\a)$ using the perfect sampling tools developed in the previous section.

\subsection{Theoretical study}

\subsubsection{A related model: the ``flip-if-not-all-equal'' PCA}

Let us define as in \cite{regn08}, the PCA FINAE($\alpha$) of neighborhood $V=\{-1,0,1\}$ and 
transition function $g:\{0,1\}^V\rightarrow \M(\{0,1\})$ given by
\[
g(x,y,z) = \a \delta_{\mbox{flip-if-not-all-equal}(x,y,z)} + (1-\a)\delta_y \:,
\]
where the function flip-if-not-all-equal (FINAE), corresponding to the elementary cellular automaton 178, 
is defined by 
\[
\mbox{flip-if-not-all-equal}(x,y,z) = \begin{cases} y & \mbox{if } x=y=z \\
1-y & \mbox{otherwise.}\end{cases}\:
\]
Clearly, $\delta_{0^\Z}$ and $\delta_{1^\Z}$ are invariant measures of
the PCA. 

Let us define $\mbox{flip-odd}: \{0,1\}^\Z\rightarrow \{0,1\}^\Z$ and 
$\mbox{flip-even}: \{0,1\}^\Z\rightarrow \{0,1\}^\Z$ by, for
        $x=(x_i)_{i\in\Z}$, 
\[
\mbox{flip-odd}(x)_i = \begin{cases} x_i & \mbox{if } i \mbox{ is
    even} \\
1-x_i & \mbox{if } i \mbox{ is odd} \end{cases}, \qquad
\mbox{flip-even}(x)_i = \begin{cases} 1-x_i & \mbox{if } i \mbox{ is
    even} \\
x_i & \mbox{if } i \mbox{ is odd} \end{cases} \:.
\]
If we extend flip-odd and flip-even to mappings on  $\M(\{0,1\}^\Z)$, we
  have
\[
\mbox{Majority}(\a) = \mbox{flip-odd} \circ \mbox{FINAE}(\a) \circ
\mbox{flip-even} \:.
\]
This equality can be checked on the local functions of the PCA Majority$(\alpha)$ and FINAE$(\alpha)$. One thus obtains that if $\pi$ is an invariant measure for FINAE$(\alpha)$, then 
\[
(\mbox{flip-odd}(\pi)+\mbox{flip-even}(\pi))/2
\]
is an invariant measure for Majority$(\alpha)$. The invariant measures $\delta_{0^\Z}$
and $\delta_{1^\Z}$ of FINAE$(\alpha)$ correspond to the invariant
measure $\mu$ in \eqref{eq-trivial} for Majority$(\alpha)$, and the existence of a
non-trivial invariant measure for FINAE$(\alpha)$ corresponds to the
existence of a second invariant measure for Majority$(\alpha)$. 

\subsubsection{Validity of the conjecture for large values of $\alpha$}

The partial result of Proposition \ref{pr-perco} relies on ideas
 from Regnault~\cite{regn08}. 

\begin{prop}\label{pr-perco} 
Let $p_c$ be the percolation threshold of directed bond-percolation in $\N^2$. 
If $\a\geq \sqrt[3]{1-(1-p_c)^4}$, then Majority$(\a)$ has several invariant measures (resp. FINAE$(\a)$ has other invariant measures than the combinations of $\delta_{0^\Z}$ and $\delta_{1^\Z}$). It is in particular the case if $\a\geq 0.996$.
\end{prop} 

{\it Proof.}
It is known that $0.6298\leq p_c\leq 2/3$, see for instance
Grimmett~\cite{grim99}. This provides the bound $\sqrt[3]{1-(1-p_c)^4}
\leq 0.996$. 


Let us consider the directed graph $G=(N,A)$ such that the set of nodes is $N=2\Z\times 2\N \cup (2\Z+1)\times(2\N+1)$
 and for each $(i,j)\in N$, there is an arc (oriented bond) from $(i,j)$ to $(i-1,j+1)$ and one from $(i,j)$ to $(i+1,j+1)$. 

\begin{figure}
\begin{center}
\includegraphics[scale=0.5]{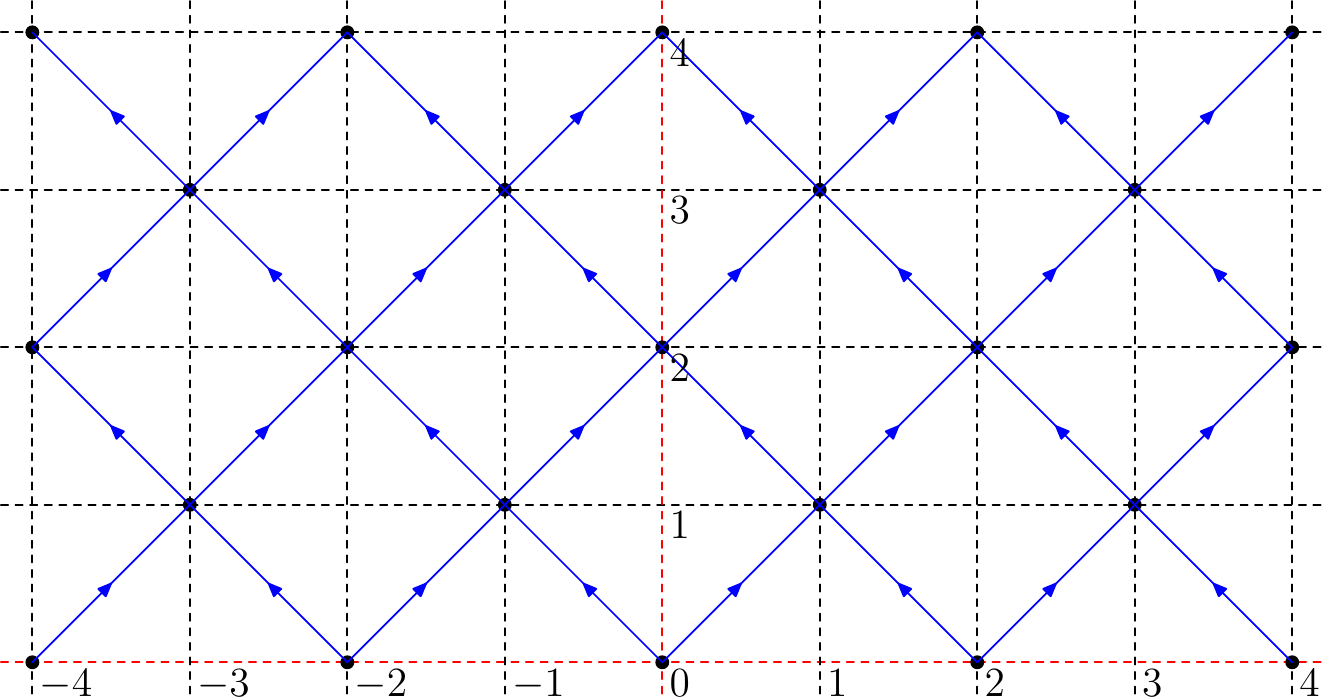} 
\caption{The graph $G$}
\end{center}
\end{figure}
 
 Let $S$ be some subset of $2\Z\times\{0\}$ called the
 \emph{source}. The oriented bond-percolation on $G$ of parameter $p$
 and source $S$ is defined as follows: each  node is open with probability $p$
 and closed with probability $1-p$, independently of the others, and a
 node of $N$ is said to be \emph{wet} if there is an open path joining
 it from some node of $S$. 

We say that the space-time diagram
 $(x^t_k)_{k\in\Z, t\in\N}$ of FINAE$(\a)$ and the percolation model 
satisfy the \emph{correspondence criterion at time
   $t$} if for each wet cell $(k,t)$ of height $t$, we
 have $x_k^ t\not=x_{k+1}^ t$ or $x_k^ t\not=x_{k-1}^ t$.    

 For values of $(\a,p)$ satisfying $\alpha\geq\sqrt[3]{1-(1-p)^4}$, Regnault is able to construct a coupling between FINAE$(\a)$ and the percolation model such that if the correspondence criterion is true
at time $t$, it is still true at time $t+1$. 
Let us take for the initial configuration of FINAE$(\a)$ the
configuration $x^0$ defined by $x^0_k=1$ if $k$ is odd and $x^0_k=0$ if $n$
is even. We also choose $S=2\Z\times\{0\}$ for the percolation model.
The correspondence criterion is true at time $0$. By the coupling described in \cite{regn08}, the criterion
is true at all time. 

Consider the percolation model and the probability $\P((0,2t) \mbox{ is wet})$. It is known (see for example \cite{grim99})
that if $p$ is strictly greater than a  certain critical value $p_c$, this probability, which decreases with $t$,
does not tend to $0$. Thus, for $p>p_c$, there exists $\eta_p>0$ such
that $\P((0,2t) \mbox{ is wet})>\eta_p$ for all $t\in\N$. By construction of the coupling,
 we obtain $\P(x_0^{2t}\not=x_1^{2t} \mbox{ or } x_0^{2t}\not=x_{-1}^{2t})\geq \eta_p$ for all
$t\in\N$.  
This proves that for $\a\geq \sqrt[3]{1-(1-p_c)^4}$, the PCA FINAE$(\a)$ has at least one invariant measure 
which is not in the convex hull of the
Dirac masses at the configurations ``all zeroes'' and ``all
ones'' 
(take any accumulation point of the Ces\`aro sums obtained from the sequence obtained from the iterated of $\delta_{x^0}$ by FINAE).
%
This result can be translated to the PCA Majority. $\hfill \square$

\subsubsection{A duality result with the double branching annihilating random walk}

The aim of this subsection is to prove a duality result between FINAE$(\a)$ and a double branching annihilating random walk (DBARW). The connection between these two models is interesting in itself and could provide a new way to study the PCA Majority$(\a)$ for small values of $\a$. A similar duality result was already obtained for interacting particle systems (see \cite{CoxDurrett}), and the behavior of the DBARW is very well understood in continuous time (see \cite{Sudbury}), but its study appears to be more difficult in discrete time.

\medskip

We now assume that $\alpha\leq 2/3$ (in particular, Proposition \ref{pr-perco} does not apply).

Let us define a process $(x_k^t)_{k\in\Z,t\in\N}$ in the following way.
For each $(k,t)\in \Z\times \N$, we first choose independently to do one (and only one) of the following:
\begin{enumerate}
\item with probability $\a/2$, draw an arc from $(k-1,t)$ to $(k,t+1)$,
\item with probability $\a/2$, draw an arc from $(k+1,t)$ to $(k,t+1)$,
\item with probability $\a/2$, draw an arc from $(k-1,t)$ to $(k,t+1)$, an arc from $(k,t)$ to $(k,t+1)$, and an arc from $(k+1,t)$ to $(k,t+1)$,
\item with probability $1-3\a/2$, draw an arc from $(k,t)$ to $(k,t+1)$.
\end{enumerate}

\begin{figure}[H]
\begin{center}
\includegraphics[]{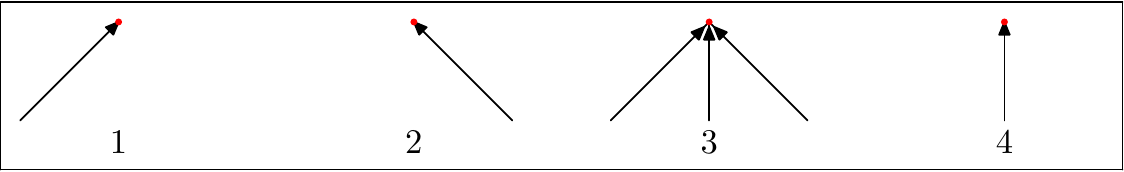}
\caption{Construction of the graph $G_1$}
\end{center}
\end{figure}

We thus obtain a directed graph $G_1$, that we will use to label each node of $\Z\times \N$  with a letter of $\{0,1\}$. 
The nodes of $\Z\times \{0\}$ are labeled according to the wanted initial configuration $x^0$. A node labeled by a $1$ will be interpreted as being occupied. A node $(k,t)\in\Z\times \N$ is then labeled by a $1$ if and only if there is an odd number of paths leading to this node from an occupied node of $\Z\times \{0\}$. This define a random field $(x_k^ t)_{k\in\Z,t\in\N}$ representing the labels of the nodes. 

\medskip

We claim that this field has the same distribution as a space-time diagram of FINAE$(\a)$ starting from $x^0$. Indeed, the value $x_k^{t+1}$ is equal to $x_{k-1}^t$ with probability $\a/2$, to $x_{k+1}^t$ with probability $\a/2$, to $x_{k-1}^t+x_{k}^t+x_{k+1}^t \mod 2$ with probability $\a/2$ and to $x_{k}^t$ with probability $1-3\a/2$. And one can check for each value $(x,y,z)\in\{0,1\}^3$ that these probabilities coincide with the ones obtained with the local function flip-if-not-all-equal. For example, if $(x_{k-1}^t,x_{k}^t,x_{k+1}^t)=(0,0,1)$, the value of $x_k^{t+1}$ will be $1$ if and only if case $2$ or case $3$ occurs, and they have together a probability $\a/2+\a/2=\a$. If $(x_{k-1}^t,x_{k}^t,x_{k+1}^t)=(0,1,0)$, we will have $x_k^{t+1}=1$ if and only if case $3$ or case $4$ occurs, which has a probability $\a/2+(1-3\a/2)=1-\a$. And if $(x_{k-1}^t,x_{k}^t,x_{k+1}^t)=(0,0,0)$ (resp. $(1,1,1)$), we will get a $0$ (resp. a $1$) in all cases.

\medskip

We now consider the process $(y_k^t)_{k\in\Z,t\in\N}$ obtained from $(x_k^t)_{k\in\Z,t\in\N}$ by reversing time. Formally, for each $(k,t)\in \Z\times \N$, we first choose independently to do one (and only one) of the following things:
\begin{enumerate}
\item with probability $\a/2$, draw an arc from $(k,t)$ to $(k-1,t+1)$;
\item with probability $\a/2$, draw an arc from $(k,t)$ to $(k+1,t+1)$;
\item with probability $\a/2$, draw an arc from $(k,t)$ to $(k-1,t+1)$, an arc from $(k,t)$ to $(k,t+1)$, and an arc from $(k,t)$ to $(k+1,t+1)$;
\item with probability $1-3\a/2$, draw an arc from $(k,t)$ to $(k,t+1)$.
\end{enumerate}

\begin{figure}[H]
\begin{center}
\includegraphics[]{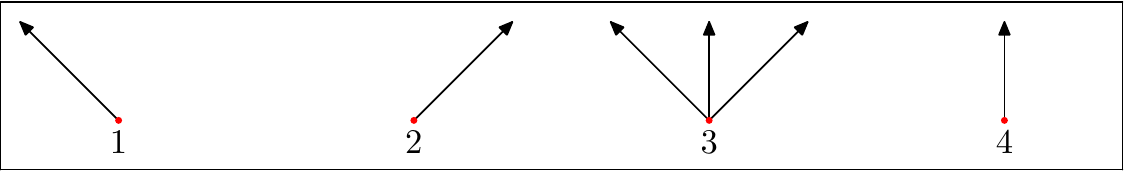}
\caption{Construction of the graph $G_2$}
\end{center}
\end{figure}

We thus obtain again a directed graph $G_2$, that we will use to label each node of $\Z\times \N$  with a letter of $\{0,1\}$. 
The nodes of $\Z\times \{0\}$ are labeled according to the wanted initial configuration $y^0$. A node labeled by a $1$ will be interpreted as being occupied. A node $(k,t)\in\Z\times \N$ is then labeled by a $1$ if and only if there is an odd number of paths leading to this node from an occupied node of $\Z\times \{0\}$. This define a random field $(y_k^ t)_{k\in\Z,t\in\N}$ representing the labels of the nodes. 

\medskip

We claim that this field has the same distribution as the double branching annihilating random walk c that we now define.

At time $0$, a particles is placed on each cell $k$ of $\Z$ such that $y^0_k=1$, and at each step of time, every particle chooses independently of the others do one (and only one) of the following things:
\begin{enumerate}
\item with probability $\a/2$, move from node $k$ to $k-1$;
\item with probability $\a/2$, move from node $k$ to $k+1$;
\item with probability $\a/2$, stay at node $k$ and create two new particles at nodes $k-1$ and $k+1$;
\item with probability $1-3\a/2$, stay at node $k$.
\end{enumerate}
If after these choices, there is an even number of particles at a node, then all these particles annihilate. If there is an odd number of them, only one particle survives. We set $w_k^t=1$ if and only if at time $t$, there is a particle at node $k$. 

To summarize, we have the following relations:
$$\mbox{ FINAE } \sim (x_k^t)_{k\in\Z,t\in\N} \stackrel{\mbox{ time-reversal }}{\longleftrightarrow}  (y_k^t)_{k\in\Z,t\in\N} \sim \mbox{ DBARW}$$

The processes $(x_k^t)_{k\in\Z,t\in\N}$ and $(y^t_k)_{k\in\Z,t\in\N}$ are obtained one from another by reversing time. This can be used to get nontrivial information for FINAE. For instance, if $A$ represents the set of occupied nodes at time $0$ for $x$, that is to say $x^0={1}_A$, we have the following duality relation:
$$\P^{\; x^0={1}_A}(x_k^t\not =x_l^t) = \P^{\; x^0={1}_A}(x_k^t+x_l^t=1)$$
$$=\P(\mbox{the total number of paths in } G_1 \mbox{ leading from } A\times\{0\} \mbox{ to } (k,t) \mbox{ or } (l,t) \mbox{ is odd})$$
$$=\P(\mbox{the total number of paths in } G_2 \mbox{ leading from } (k,0) \mbox{ or } (l,0) \mbox{ to } A\times\{t\} \mbox{ is odd})$$
$$=\P^{\; y^0={1}_{\{k,l\}}}(\sum_{i\in A} y^t_i \mbox{ is odd})$$
$$\leq \P^{\; y^0={1}_{\{k,l\}}}(\exists i\in A, y^ t_i=1).$$

Thus, to prove that the probability for the PCA FINAE that two cells
$k$ and $l$ will be in different states at time $t$ tends to $0$  as
$t$ tends to $+\infty$, it is sufficient to prove that in the DBARW,
starting from two particles, the probability of extinction of the
population of particles tends to $1$. 

\subsection{Experimental study}

We tried to get some numerical evidence for Conjecture \ref{co-co}. To study the PCA Majority experimentally, a first idea would be to consider the
same PCA on the set of cells $\Z_n$, $n$ odd. This does not work
well. 
First, due to the state space explosion, computing exactly the invariant
measure is possible only for small values (we did it up to $n=9$ using Maple). 
Second, the algorithms of Section \ref{sec:sampling} cannot be applied since the EPCA is
not ergodic. 

Instead, we use approximations of the PCA by NH-PCA on a
finite subset of cells, the methodology sketched in Section
\ref{sse-exte}. Again, computing exactly the invariant measure is
impossible except for very small windows. But now the sampling
algorithms become effective. 

\medskip

\begin{figure}
\begin{multicols}{2}
\begin{center}
\vfill~

\subfloat[The value of $c_n$ as a function of  $n$, for different $\alpha$.]{\includegraphics[trim = 2cm 1.5cm 2cm 2cm, clip, width=9cm]{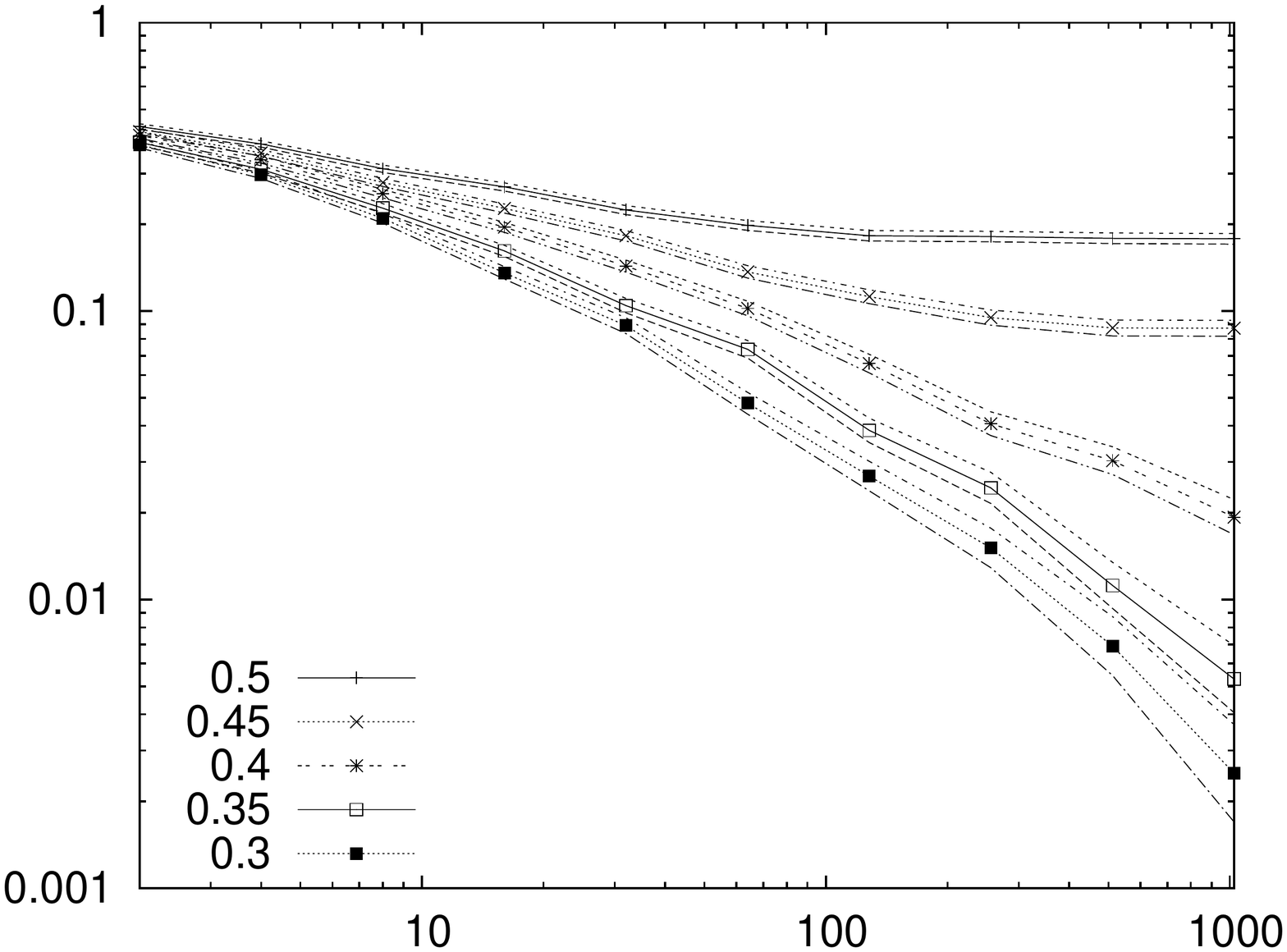}}

\vfill~\subfloat[$\alpha=0.5$]{\includegraphics[trim = 10cm 0.1cm 10cm 0.1cm,
    clip, width=5cm]{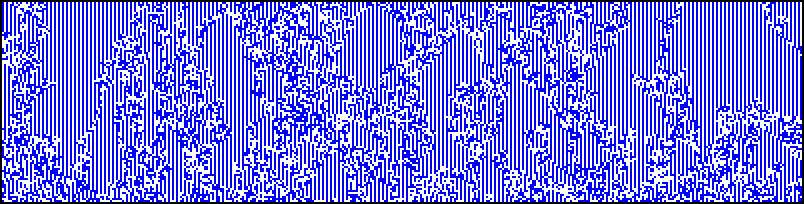}} 

\vfill~\subfloat[$\alpha=0.3$]{\includegraphics[trim = 10cm 0.1cm 10cm 0.1cm,
    clip, width=5cm]{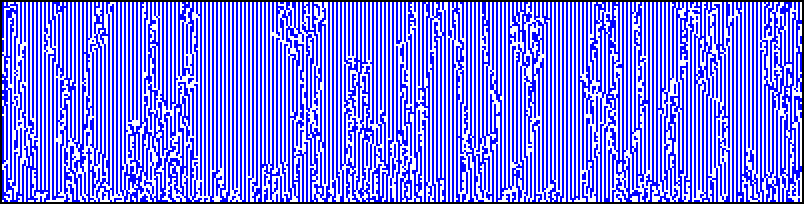}} 
\end{center}
\end{multicols}
\caption{Experimental study of Majority($\alpha)$ (the configurations
  at odd times only are represented on the space-time
  diagrams). \label{fi-exper}} 
\end{figure}

Let $P$ be the PCA Majority. 
Set $D_n=\{-n,\ldots,n\}$, and let $\nu$ be the uniform measure on
$\{0,1\}$. Consider the NH-PCA $P(\nu,D_n)$. Let $\mu_n$ be the unique
invariant measure of $P(\nu,D_n)$.  We are interested in the quantity
\[
c_n=\mu_n\{x\in X \mid
x_0=x_1=0\} +  \mu_n \{x\in X \mid x_0=x_1=1\} \:.
\]
Indeed, by application of
Lemma \ref{pr-approx}, if $\limsup_n c_n >0$, then there exists a
non-trivial invariant measure for the PCA Majority on $\Z$. 



Now the NH-EPCA is ergodic, so the sampling algorithms of Section \ref{sec:sampling} can
be used. We were able to run the algorithms up to a window size of $n=1024$
before running into overtime problems. 
The experimental results appear in Figure \ref{fi-exper}, with a logarithmic scale. 
We ran the sampling algorithms $10 000$ times.  We show on the figure the confidence intervals
calculated with Wilson score test at 95\%. 

\medskip

It is reasonable to believe that the top
two curves in Figure \ref{fi-exper} do not converge to 0 while the
bottom three converge to 0.  
This is consistent with the visual impression
of space-time diagrams. It reinforces Conjecture \ref{co-co} with
a possible phase transition between 0.4 and 0.45. 

\paragraph{Acknowledgements.}

We used the applet FiatLux developed by N. Fat\`es and available on
his website (LORIA, INRIA Lorraine) to draw the space-time diagrams of
Figures \ref{fi-pca1} and \ref{fi-exper}.

\FloatBarrier

\addcontentsline{toc}{section}{References}


\end{document}